\documentclass{amsart}
\usepackage{amscd,amssymb,stmaryrd}
\usepackage[all]{xy}
\renewcommand{\mod}{\operatorname{mod}\nolimits}
\newcommand{\Mod}{\operatorname{Mod}\nolimits}
\newcommand{\ind}{\operatorname{ind}\nolimits}

\newcommand{\gr}{{\operatorname{gr}\nolimits}}
\newcommand{\tr}{{\operatorname{tr}\nolimits}}
\newcommand{\Gr}{{\operatorname{Gr}\nolimits}}

\newcommand{\add}{\operatorname{add}\nolimits}
\newcommand{\soc}{\operatorname{soc}\nolimits}

\newcommand{\Hom}{\operatorname{Hom}\nolimits}
\newcommand{\End}{\operatorname{End}\nolimits}

\renewcommand{\Im}{\operatorname{Im}\nolimits}
\newcommand{\Ker}{\operatorname{Ker}\nolimits}

\newcommand{\rrad}{\mathfrak{r}}
\newcommand{\rad}{\operatorname{rad}\nolimits}

\newcommand{\gldim}{\operatorname{gldim}\nolimits}
\newcommand{\LL}{\operatorname{LL}\nolimits}
\newcommand{\Tr}{\operatorname{Tr}\nolimits}
\newcommand{\Ext}{\operatorname{Ext}\nolimits}
\newcommand{\Tor}{\operatorname{Tor}\nolimits}

\newcommand{\op}{{\operatorname{op}\nolimits}}

\newcommand{\pd}{{\operatorname{pd}\nolimits}}

\newcommand{\fd}{{\operatorname{fd}\nolimits}}
\newcommand{\comp}{\operatorname{\scriptstyle\circ}}
\newcommand{\m}{\mathfrak{m}}

\newcommand{\frakm}{\mathfrak{m}}

\newcommand{\G}{\Gamma}
\renewcommand{\L}{\Lambda}

\newcommand{\A}{{\mathcal A}}
\newcommand{\C}{{\mathcal C}}
\newcommand{\D}{{\mathcal D}}
\newcommand{\E}{{\mathcal E}}
\newcommand{\F}{{\mathcal F}}
\newcommand{\calH}{{\mathcal H}}

\renewcommand{\S}{{\mathcal S}}
\newcommand{\T}{{\mathcal T}}

\newcommand{\extto}{\xrightarrow}

\newtheorem{lem}{Lemma}[section]
\newtheorem{prop}[lem]{Proposition}
\newtheorem{cor}[lem]{Corollary}
\newtheorem{thm}[lem]{Theorem}
\theoremstyle{definition}
\newtheorem{defin}[lem]{Definition}

\newtheorem{example}[lem]{Example}

\begin{document}

\title{Artin-Schelter regular algebras and categories}
\author[Martin\'ez-Villa]{Roberto Martin\'ez-Villa}
\thanks{The first author thanks the Universidad Nacional Autonoma de
  Mexico program PAPITT for funding
  the research project. In addition he thanks his coauthor and
  Department of Mathematical Sciences (NTNU) for their kind
  hospitality and support through the NFR Storforsk grant no.\
  167130.} 
\address{Roberto Martin\'ez-Villa, Instituto de Matem\'aticas,
Universidad Nacional Autonoma de Mexico, Campus Morelia, Apartado
Postal 27-3 (Xangari), C.P. 58089, Morelia, Michoac\'an, Mexico}
\email{mvilla@matmor.unam.mx}
\author[Solberg]{\O yvind Solberg}
\thanks{The second author thanks Universidad Nacional Autonoma de
Mexico and his coauthor for their kind hospitality and support, and
the Department of Mathematical Sciences (NTNU) and the NFR Storforsk
grant no.\ 167130 also for the support.} 
\address{\O yvind Solberg\\
Institutt for matematiske fag\\
NTNU\\
N--7034 Trondheim\\
Norway}
\email{oyvinso@math.ntnu.no}
\thanks{Some of the results in this paper were presented at the
  workshop \emph{Representation Theory of Finite-Dimensional Algebras}
  held at Mathematisches Forschungsinstitut Oberwolfach, February 2005.}
\date{\today}

\begin{abstract}
  Motivated by constructions in the representation theory of finite
  dimensional algebras we generalize the notion of Artin-Schelter
  regular algebras of dimension $n$ to algebras and categories to
  include Auslander algebras and a graded analogue for infinite
  representation type. A generalized Artin-Schelter regular algebra or
  a category of dimension $n$ is shown to have common properties with
  the classical Artin-Schelter regular algebras. In particular, when
  they admit a duality, then they satisfy Serre duality formulas and
  the $\Ext$-category of nice sets of simple objects of maximal
  projective dimension $n$ is a finite length Frobenius category.
\end{abstract}

\maketitle

\section*{Introduction}

Artin-Schelter regular rings have been introduced as non-commutative
analogues of polynomial rings by Artin and Schelter in
\cite{ArS}. This has become the starting point of a rich theory of
non-commutative algebra and algebraic geometry. Our initial source of
inspiration and examples comes from the representation theory of
Artin/finite dimensional algebras. A finite dimensional algebra $\L$
has either a finite or an infinite number of isomorphism classes of
indecomposable finitely generated modules, socalled finite or infinite
type, respectively. In either case, one associates the Auslander
algebra/category $\mathfrak{A}_\L$ to the algebra $\L$, and we have
observed that $\mathfrak{A}_\L$ behaves very much like an
Artin-Schelter regular algebra, namely
\begin{enumerate}
\item[(i)] the global dimension is $2$ (finite) (\cite{A}),
\item[(ii)] all the simples of maximal projective dimension satisfy
  the Gorenstein condition, or equivalently the $2$-simple condition
  (\cite{I3}),
\item[(iii)]\sloppy the $\Ext$-algebra/category of (nice) sets of
  simples of maximal projective dimension is a finite
  dimensional/length Frobenius algebra/category.
\end{enumerate}
This motivated us to investigate this in further detail and
generality, and the purpose of this paper is to present a
generalization of Artin-Schelter regular algebras to algebras and
categories that include the Auslander algebras/categories occurring in
representation theory of Artin/finite dimensional
algebras. Furthermore, as pointed out to us by Osamu Iyama, the
additive closure of a $n$-cluster tilting module over a finite
dimensional algebra is a generalized Artin-Schelter regular category.
For rings the same generalization or class of rings has been
considered by Iyama in \cite{I3}, where Noetherian rings having finite
global dimension $n$ satisfying the $n$-simple condition have
characterized in terms of properties of the minimal injective
resolution of the ring and the opposite ring. Here the main motivation
was to study orders over complete discrete valuation rings.

No assumption on the Gelfand-Kirillov dimension is needed in Iyama's
work. Also in our definition of a generalized Artin-Schelter regular
algebra or category, we impose no requirement on having finite
Gelfand-Kirillov dimension. From our point of view this do not enter
or effect the definition, and in addition we construct examples of
generalized Artin-Schelter regular algebras which have infinite
Gelfand-Kirillov dimension. But an algebra in the class of examples we
construct has finite Gelfand-Kirillov dimension if and only if it is
Noetherian. This last question is, to our knowledge, an open problem
for Artin-Schelter regular algebras.

In this paper we freely use the results from \cite{MVS1}. Furthermore,
in a forthcoming paper we use the results in this paper to show that
the some properties of the associated graded Auslander category of a
component in the Auslander-Reiten quiver of a finite dimensional
algebra is reflected as Noetherianity and different Gelfand-Kirillov
dimensions.

Recall that a graded connected $K$-algebra $\L$ over a field $K$ is
\emph{Artin-Schelter regular of dimension $n$}, if
$\gldim\L=n<\infty$, it has finite Gelfand-Kirillov dimension, and the
unique simple graded $\L$-module $S$ satisfies the Gorenstein
condition, that is, $\Ext^n_\L(S,\L)=(0)$ for $i\neq n$ and
$\Ext^n_\L(S,\L)$ is isomorphic to some graded shift of the unique
simple graded $\L^\op$-module. As mentioned above our goal is to
extend this notion of Artin-Schelter regular algebras such that it
includes the (graded) Auslander algebras and categories. We do this in
two steps. In Section \ref{section:1} we first consider the class of
rings considered by Iyama in \cite{I3} (even semiperfect) and in the
graded setting non-connected algebras. Properties of generalized
Artin-Schelter regular algebras are discussed in Section
\ref{section1.5}, and in particular we show that the $\Ext$-algebra of
a finite set of simple modules of maximal projective dimension
permuted the functor $D\Ext^n_\L(-,\L)$ is a finite dimensional
Frobenius algebra. Section \ref{section2} is devoted to reviewing the
polynomial ring $K[x_1,x_2,\ldots,x_n]$ as an example of ungraded
generalized Artin-Schelter algebras to illustrate our results. We
construct in Section \ref{section3} a family of generalized
Artin-Schelter regular algebras including some non-Noetherian and with
infinite Gelfand-Kirillov dimension. This is contrary what is believed
to be true for classical Artin-Schelter regular algebras, but we show
that they are Noetherian if and only if they have finite
Gelfand-Kirillov dimension. We take the second step in Section
\ref{section4} and give the definition of generalized Artin-Schelter
regular $K$-categories. In Section \ref{section:6}, we show that our
generalized Artin-Schelter regular categories share similar properties
with the generalized Artin-Schelter regular algebras. The final
section, Section \ref{section:7}, is devoted to studying generalized
Artin-Schelter regular categories $\C$ of global dimension $2$, where
we show that $\C$ is coherent and $\Gr(\C)$ is abelian when $\C$ is
generated in degrees $0$ and $1$.

\section{Generalized Artin-Schelter regular algebras}\label{section:1}

This section is devoted to generalizing the classical notion of
Artin-Schelter regular algebras of Artin and Schelter from
\cite{ArS}. For ungraded rings the same class of Noetherian rings has
been considered by Iyama in \cite{I3}.  We have seen the need of such
a generalization through introducing a Koszul theory associated to any
finite dimensional algebra in \cite{MVS1} (and applications in
forthcoming papers). The results in this section serve as a model for
our generalization to categories in Section \ref{section4}.

We start the section by reviewing how the Koszul theory/object we
associate to any finite dimensional algebra in \cite{MVS1} is similar
in behavior to Artin-Schelter regular algebras. Then we develop the
necessary theory to end the section with our definition of a
generalized Artin-Schelter regular algebra.

Before indulging in the gruesome details of our generalization of
Artin-Schelter regular algebras, we first try to motivate the need for
such a generalization. For another motivation we refer the reader to
the papers \cite{I1,I2,I3}. First we consider the class of Auslander
algebras. They correspond to the finite dimensional algebras (Artin
algebras) with only a finite number of isomorphism classes of
indecomposable finitely generated left modules (finite representation
type). They are defined as finite dimensional algebras with dominant
dimension at least $2$ and global dimension at most $2$. Let $\L$ be a
finite dimensional $K$-algebra of finite representation type, and
denote by $M$ the direct sum of all indecomposable finitely
generated left modules with one from each isomorphism class. Then
$\G=\End_\L(M)^\op$ is the corresponding Auslander algebra. The simple
$\G$-modules are in one-to-one correspondence with the indecomposable
$\L$-modules. Let $\{S_C\}_{C\in\ind\L}$ be the simple $\G$-modules,
where $\ind\L$ is the set of isomorphism classes of indecomposable
finitely generated left $\L$-modules. A minimal projective resolution
of of $S_C$ is given by
\[0\to \Hom_\L(M,\tau C)\to \Hom_\L(M,B)\to \Hom_\L(M,C)\to S_C\to 0\]
if $C$ is non-projective and $0\to \tau C\to B\to C\to 0$ is the
almost split sequence ending in $C$, and by
\[0\to \Hom_\L(M,\rrad P)\to \Hom_\L(M,P)\to S_P\to 0\] if $C=P$ is
projective. Here $\rrad$ denotes the Jacobson radical of $\L$. For the
definition and further properties of almost split 
sequences we refer to \cite{ARS}. It is easy to see that all the
simple modules corresponding to non-projective indecomposable modules
satisfy the Gorenstein condition or the $2$-simple condition, that is,
$\Ext^i_\G(S_C,\G)=(0)$ for $i\neq 2$ and $\Ext^2_\G(S_C,\G)$ is a
simple $\G^\op$-module. Note also that they have the highest possible
projective dimension $2$, as the global dimension of $\Gamma$ is
$2$. Hence, in this case, all the simple modules of highest possible
projective dimension behave as over an Artin-Schelter regular
algebra. This is called the $n$-simple condition in
\cite{I3}. Similarly one can consider the replacement of the Auslander
algebra for a finite dimensional algebra of infinite representation
type, namely, the category of additive functors from $(\mod\L)^\op$ to
$\Mod K$. Here one can find similar behavior. In these situations we
are able to prove similar results as for the classical Artin-Schelter
regular algebras. This motivates our generalization of Artin-Schelter
regular algebras and also Artin-Schelter regular $K$-categories.

To exemplify what we mean by similar results we give a concrete
example.

\begin{example}\label{exam:Auslanderalgebra}
Let $\Sigma=K[x]/(x^3)$ for a field $K$. Then $\Sigma$ has three
indecomposable finitely generated modules, $U=K[x]/(x)$,
$L=K[x]/(x^2)$ and $\Sigma$. Let $M=U\amalg L\amalg \Sigma$, then
$\L=\End_\Sigma(M)^\op$ is the corresponding Auslander algebra for
$\Sigma$. The algebra $\L$ is a finite dimensional algebra over $K$
with three simple modules $S_U$, $S_L$ and $S_\Sigma$. The exact
sequences $0\to U\to L\to U\to 0$, $0\to L\to \Sigma\amalg U\to L\to
0$ and $0\to L\to \Sigma$ induce the following projective resolutions
\[0\to \Hom_\Sigma(M,U)\to \Hom_\Sigma(M,L)\to \Hom_\Sigma(M,U)\to S_U\to 0,\]
\[0\to \Hom_\Sigma(M,L)\to \Hom_\Sigma(M,\Sigma\amalg U)\to
\Hom_\Sigma(M,L)\to S_L\to 0,\]
and
\[0\to \Hom_\Sigma(M,L)\to \Hom_\Sigma(M,\Sigma)\to S_\Sigma\to 0\] of
the three simple $\L$-modules $S_U$, $S_L$ and $S_\Sigma$.  By
applying the functor $\Hom_\L(-,\L)$ to the first two sequences, it
follows easily that $\Ext_\L^i(S,\L)=(0)$ for $i\neq 2$, and that
$\Ext_\L^2(S,\L)\simeq S^\op$ for $S$ equal to $S_U$ or $S_L$. Hence
the simples $S_U$ and $S_L$ satisfies the $2$-simple condition, while
$S_\Sigma$ does not. Now let $T=S_U\amalg S_L$, and consider
$\G=\Ext^*_\L(T,T)$. It is easy to see that $\G$ is isomorphic to the
path algebra $K(\xymatrix{1\ar@<1ex>[r]^\alpha &
  2\ar@<1ex>[l]^\beta})/(\alpha\beta\alpha, \beta\alpha\beta)$, which
is Frobenius. Recall that for an Artin-Schelter regular algebra the
$\Ext$-algebra of the simple graded module is Frobenius.  Hence the
above example illustrate a similar behavior, and for our generalized
Artin-Schelter regular algebras of dimension $n$ we show that the
$\Ext$-algebra of certain subsets of simple modules satisfying the
$n$-simple condition (also over the opposite algebra) is finite
dimensional Frobenius.
\end{example}

Now we start discussing the generalization of the notion of
Artin-Schelter regular algebras.  Let $\L$ be any ring and $n$ a
positive integer. We denote by $\calH^n_\L$ the full subcategory of
all left $\L$-modules $M$ with projective dimension $n$ such that
$\Ext^i_\L(M,\L)=(0)$ for $i\neq n$ and with a projective resolution
consisting of finitely generated $\L$-modules. Similar subcategories
of $\L$-modules are considered in \cite{I3}, and similar results as
below are proven. For completeness we include the proofs here. Let
$\tr^n_\L=\Ext^n_\L(-,\L)\colon \Mod\L\to \Mod\L^\op$. Denote by
$(-)^*$ the functor $\Hom_\L(-,\L)\colon \Mod\L\to \Mod\L^\op$.
\begin{prop}\label{prop:duality}
\begin{enumerate}
\item[(a)] \sloppy The subcategory $\calH_\L^n$ is closed under
extensions. If $\gldim\L = n$, then $\calH_\L^n$ is closed under
cokernels of monomorphisms.
\item[(b)] The functor $\tr^n_\L\colon \Mod\L\to \Mod\L^\op$ restricts
to a functor $\tr^n_\L\colon \calH^n_\L\to \calH^n_{\L^\op}$, and here
it is an exact duality.
\item[(c)] For any $\L$-module $M$ in $\calH^n$ and any $\L$-module
$L$ there are isomorphisms
\[\Tor^\L_i(\tr^n_\L(M),L) \simeq \Ext^{n-i}_\L(M,L)\]
for all $i\geq 0$.
\end{enumerate}
\end{prop}
\begin{proof}
Let $\eta\colon 0\to A\to B\to C\to 0$ be an exact sequence in
$\Mod\L$.

(a) If $A$ and $C$ are in $\calH^n_\L$, then the long exact sequence
induced by $\eta$ and the Horseshoe Lemma imply that $B$ is in
$\calH^n_\L$. Hence $\calH^n_\L$ is closed under extensions.

Suppose that $\gldim\L=n$. In the exact sequence $\eta$ assume that
$A$ and $B$ are in $\calH^n_\L$. Then $C$ also has finite projective
dimension. Using the Horseshoe Lemma on the exact sequence $0\to
\Omega^1_\L(B)\to \Omega^1_\L(C)\amalg P\to A\to 0$ for some
projective $P$, imply that $C$ has a finitely generated projective
resolution.  The long exact sequence induced by $\eta$ then imply that
$C$ is in $\calH^n_\L$.

(b) Let $M$ be in $\calH^n_\L$. Let
\[0\to P_n\to P_{n-1}\to \cdots \to P_1\to P_0\to M\to 0\]
be a finitely generated projective resolution of $M$. Since $M$ is in
$\calH^n_\L$, the above long exact sequence gives rise to the long
exact sequence
\[0\to P^*_0\to P^*_1\to \cdots \to P^*_{n-1}\to P^*_n\to
\Ext^n_\L(M,\L)\to 0\]
This shows that $\tr^n_\L(M)$ has a finitely generated projective
resolution of length $n$. Applying the functor $\Hom_{\L^\op}(-,\L)$,
which we also denote by $(-)^*$, we obtain the following commutative
diagram
\[\xymatrix@C=10pt{
& 0\ar[r] & P_n\ar[r]\ar[d]
        & P_{n-1}\ar[r]\ar[d]
        & \cdots\ar[r]
        & P_1\ar[r]\ar[d]
        & P_0\ar[r]\ar[d]
        & M\ar[r]\ar[d]
        & 0\\
0\ar[r]& (\tr^n_\L(M))^*\ar[r] & P_n^{**} \ar[r]
        & P_{n-1}^{**}\ar[r]
        & \cdots\ar[r]
        & P_1^{**}\ar[r]
        & P_0^{**}\ar[r]
        & \tr^n_{\L^\op}(\tr^n_\L(M)) \ar[r]
        & 0
}\]
where all vertical maps $P_i\to P_i^{**}$ are isomorphisms. In
addition, $P_1^{**}\to P_0^{**}\to \tr^n_{\L^\op}(\tr^n_\L(M))\to 0$ is
exact. It follows from this that $\Ext^i_{\L^\op}(\tr^n_\L(M),\L)=(0)$ for
all $i\neq n$. Hence $\tr^n_\L(M)$ is in $\calH^n_{\L^\op}$ and $M\simeq
\tr^n_{\L^\op}(\tr^n_\L(M))$. This shows that $\tr^n_\L\colon \calH^n_\L\to
\calH^n_{\L^\op}$ is a duality.

(c) Let $M$ be in $\calH^n_\L$ with projective resolution as above in
(b). Then the complexes $0\to P^*_0\otimes_\L L\to P^*_1\otimes_\L L\to
\cdots \to P_{n-1}^*\otimes_\L L\to P_n^*\otimes_\L L\to 0$ and $0\to
\Hom_\L(P_0,L)\to \Hom_\L(P_1,L)\to \cdots\to \Hom_\L(P_{n-1},L)\to
\Hom_\L(P_n,L)\to 0$ are isomorphic. The claims follow directly from
this.
\end{proof}

If a simple module $S$ is in $\calH^n_{\L^\op}$, it is not clear that
$\tr^n_{\L^\op}(S)$ in $\calH^n_\L$, again is a simple module. For
Noetherian rings this property is studied in \cite{I3} . There a
Noetherian ring $\L$ is said \emph{to satisfy the $n$-simple
  condition} if, for every simple $\L$-module $S$ with $\pd_\L S=n$,
we have that $\Ext^i_\L(S,\L)=(0)$ for $0\leq i < n$ and
$\Ext^n_\L(S,\L)$ is a simple $\L^\op$-module. The following follows
directly from \cite[Proposition 6.3]{I3}.
\begin{thm}[\cite{I3}]
Let $\L$ be a Noetherian ring with $\gldim \L\leq n$. Then $\L$ and
$\L^\op$ satisfy the $n$-simple condition if and only if, in a minimal
injective resolution
\[0\to \L\to I^0\to I^1\to I^2\to \cdots\]
we have $\fd_\L I^i<n$ for $0\leq i<n$ and the same property for
$\L^\op$.
\end{thm}
We do not know if there is a similar characterization of the
$n$-simple condition for non-Noetherian rings. It is shown in
\cite[The proof of Theorem 2.10]{I4} that $\Ext^n_\L(S,\L)$ has finite
length over $\L^\op$, whenever $S$ is in $\calH^n_\L$. So there is
always simple submodules of $\Ext^n_\L(S,\L)$ for $S$ in
$\calH^n_\L$. We show next that, if one of these simples are in
$\calH^n_{\L^\op}$, then this is enough for having the $n$-simple
condition.

\begin{lem}\label{lem:simpletosimple}
Let $\L$ be a ring of finite global dimension $n$, and let $S$ be a
simple module in $\calH^n_{\L^\op}$. Assume that $\tr^n_{\L^\op}(S)$
contains a simple module in $\calH^n_\L$ as a submodule. Then
$\tr^n_{\L^\op}(S)$ is a simple module (and it is in $\calH^n_\L$).
\end{lem}
\begin{proof}
Let $S$ be a simple module in $\calH^n_{\L^\op}$. Assume that $T\subseteq
\tr^n_{\L^\op}(S)$ is a simple submodule with $T$ in $\calH^n_\L$. By
Proposition \ref{prop:duality} we obtain an exact sequence
\[0\to \tr^n_\L(\tr^n_{\L^\op}(S)/T)\to \tr^n_\L(\tr^n_{\L^\op}(S))\to
\tr^n_\L(T)\to 0\]
in $\calH^n_\L$, where $\tr^n_\L(\tr^n_{\L^\op}(S))\simeq S$ and
$\tr^n_\L(T)$ is nonzero. It follows that $S\simeq \tr^n_\L(T)$ and
$\tr^n_{\L^\op}(S)\simeq T$.
\end{proof}

Let $\S^n$ denote the simple $\L$-modules which are in
$\calH^n_\L$. Then $\calH^n_\L$ contains all $\L$-modules of finite
length with composition factors only in $\S^n$ by Proposition
\ref{prop:duality}. Our next aim is to find sufficient conditions on
$\L$ such that $\calH^n_\L$ is exactly the full category of
$\L$-modules of finite length with composition factors only in $\S^n$.

Define a subfunctor $t^n\colon \Mod\L\to \Mod\L$ of the identity as
follows. Let $M$ be a $\L$-module. Then let $t^n(M)=\sum_{L\subseteq
M} L$, where the sum is taken over all submodules $L$ of $M$, where
$L$ has finite length with composition factors only in $\S^n$. We say
that a module $M$ is \emph{$\calH^n_\L$-torsion} if $M=t^n(M)$. A
module $M$ is \emph{$\calH^n_\L$-torsion free} if $t^n(M)=(0)$. Next
we show that $t^n$ is a radical in $\Mod\L$ for all $n\geq 1$.

\begin{lem}\label{lem:toftorsionfree}
Let $\L$ be a semiperfect left Noetherian ring, or let $\L$ be a
positively graded ring such that $\L_0\simeq K^t$ for some $t$,
$\dim_K\L_i<\infty$ for all $i\geq 0$, and $\L_{\geq 1}$ is a finitely
generated left $\L$-module. Then $t^n(M/t^n(M))=(0)$ for any
$\L$-module or graded $\L$-module $M$, or equivalently $M/t^n(M)$ is
$\calH^n_\L$-torsion free.
\end{lem}
\begin{proof}
Let $M$ be a $\L$-module.  Assume that $M/t^n(M)$ has a submodule of
finite length with compositions factors only in $\S^n$. Then $M/t(M)$
has a simple submodule $S$ from $\S^n$, where $S=\L/\m$ for some
maximal left (graded) ideal $\m$. Let $x$ in $M$ be such that
$\overline{x}=1+\m$. Then $\m x$ is contained in $t^n(M)$.  If $\L$ is
a left Noetherian ring, then $\m$ is finitely generated. If $\L$ is
graded as above, then $\m=\m_0+\L_{\geq 1}$ for some maximal ideal
$\m_0$ in $\L_0$. Hence, also here $\m$ is finitely
generated under the assumptions in the graded case.  In either case,
it follows that $\m x$ is contained in a submodule $L$ of $M$ of
finite length with composition factors only in $\S^n$. Furthermore,
there is an exact sequence $0\to \m x\to \L x\to S\to 0$, so that $\L
x$ is submodule of $M$ of finite length with composition factors only
in $\S^n$. Hence $x$ is in $t^n(M)$. This is a contradiction, so we
infer that $t^n(M/t^n(M))=(0)$.
\end{proof}

In some settings torsion free is the same as a submodule of a free
module. If, in addition, the ring has finite global dimension, then
the projective dimension of a torsion free module is always at most
one less than the global dimension. We show next that the behavior of
$\calH^n_\L$-torsion free modules are similar when $\L$ satisfies some
additional properties compared to Lemma \ref{lem:toftorsionfree}.

\begin{prop}\label{prop:torsionfreepddim}
Let $\L$ be a semiperfect left Noetherian ring, or let $\L$ be a
positively graded ring such that $\L_0\simeq K^t$ for some $t$,
$\dim_K\L_i<\infty$ for all $i\geq 0$, and $\L_{\geq 1}$ is a finitely
generated left $\L$-module.

Assume that $\L$ has finite (graded) global dimension $n$, that the
category $\calH^n_{\L^\op}$ contains all (graded) simple
$\L^\op$-modules of maximal projective dimension $n$ and that
$\tr^n_{\L^\op}(S)$ has a simple module from $\calH^n_\L$ as a
submodule for all (graded) simple $\L^\op$-modules $S$ of maximal
projective dimension $n$. Then any finitely generated
$\calH^n_\L$-torsion free $\L$-module has projective dimension at most
$n-1$.
\end{prop}
\begin{proof}
Let $\S^n_{\L^\op}$ be all the (graded) simple $\L^\op$-modules of
maximal projective dimension $n$, and assume that $\S^n_{\L^\op}$ is
in $\calH^n_{\L^\op}$.  Let $S'=\tr^n_{\L^\op}(S)$ for a simple
$\L^\op$-module $S$ in $\S^n_{\L^\op}$.  By Proposition
\ref{prop:duality} (c) we have that $\Ext^{n-i}_\L(S',L)\simeq
\Tor^\L_i(S,L)$ for all $i$ and for any $\L$-module $L$. In
particular, assuming that $L$ is $\calH^n_\L$-torsion free and by
Lemma \ref{lem:simpletosimple} the module $S'$ is a simple module in
$\calH^n_\L$, we obtain
\[(0)=\Hom_\L(S',L)\simeq \Tor^\L_n(S,L)\simeq
\Tor^\L_1(S,\Omega^{n-1}_\L(L))\]
for simple $\L^\op$-modules $S$ in $\S^n_{\L^\op}$. For any other
simple $\L^\op$-module $S''$ the projective dimension of $S''$ is at
most $n-1$, so that $\Tor^\L_1(S'',\Omega^{n-1}_\L(L))=(0)$ also for
these simples. Hence $\Tor^\L_1(S,\Omega^{n-1}_\L(L))=(0)$ for all
simple $\L^\op$-modules $S$.

If $\L$ is graded as in the claim of the proposition, then
$\Omega^{n-1}_\L(L)$ is a graded $\L$-module bounded below. Let $0\to
\Omega^n_\L(L)\to P\to \Omega^{n-1}(L)\to 0$ be a graded projective
cover. Then $(0)=\Tor^\L_1(\L/\rad\L,\Omega^{n-1}_\L(L))\simeq
\Omega^n_\L(L)/\rad\L \Omega^n_\L(L)$, where $\rad\L$ is the graded
Jacobson radical of $\L$. Each graded part of $\Omega^n_\L(L)$ is a
finitely generated $\L_0$-module. Using Nakayama's Lemma we infer that
$\Omega^n_\L(L)=(0)$.

If $\L$ is left Noetherian, then $\Omega^{n-1}_\L(L)$ is a finitely
generated flat $\L$-module, hence it is projective. In either case it
follows that the projective dimension of $L$ is at most $n-1$.
\end{proof}

A $\calH^n_\L$-torsion module $M$ need not to be in $\calH^n_\L$ (not
even finitely generated), but they share the property that
$\Ext^i_\L(-,\L)$ vanish for all $i\neq n$ with the modules in
$\calH^n_\L$, as is shown next.

\begin{lem}\label{lem:torsioninH}
Let $\L$ be a ring of finite global dimension $n$ such that all simple
$\L$-modules of maximal projective dimension $n$ are in
$\calH^n_\L$. Assume that $\L$ is a semiperfect left Noetherian ring,
or that $\L$ is a positively graded ring such that $\L_0\simeq K^t$
for some $t$, $\dim_K\L_i<\infty$ for all $i\geq 0$, and $\L_{\geq 1}$
is a finitely generated left $\L$-module.

Then any $\calH^n_\L$-torsion module $M$ satisfies
$\Ext^i_\L(M,\L)=(0)$ for all $i\neq n$.
\end{lem}
\begin{proof}
Let $M$ be a $\calH^n_\L$-torsion $\L$-module. Then $M$ can be written
as a filtered colimit of submodules $M_\alpha$ of finite length with
composition factors only in $\S^n$, such that $M_\alpha\hookrightarrow
M_{\alpha'}$ whenever $\alpha\leq \alpha'$.  Let $0\to \L\to I_0\to
I_1\to \cdots \to I_{n-1}\to I_n\to 0$ be a minimal injective
resolution of $\L$. Since
\[0\to (M_\alpha,\Omega^{-j}_\L(\L))\to (M_\alpha,I_j) \to
(M_\alpha,\Omega^{-j-1}_\L(\L)) \to 0\]
induces an exact sequence of inverse system where the first is a
surjective system for $j=0,1,\ldots,n-2$, the inverse limit of each of
these different systems are exact by the Mittag-Leffler
condition. Hence they give rise to the exact sequences
\[0\to (M,\Omega^{-j}_\L(\L))\to (M,I_j) \to
(M,\Omega^{-j-1}_\L(\L)) \to 0\]
for $j=0,1,\ldots,n-2$, and therefore $\Ext^j_\L(M,\L)=(0)$ for all
$j\neq n$.
\end{proof}

Now we are ready to give sufficient conditions for $\calH^n_\L$ to be
exactly the modules of finite length with composition factors only in
$\S^n$.

\begin{prop}\label{prop:exactlytorsion}
Let $\L$ be a ring of finite global dimension $n$. Assume that $\L$ is
a semiperfect left Noetherian ring, or that $\L$ is a positively
graded ring such that $\L_0\simeq K^t$ for some $t$,
$\dim_K\L_i<\infty$ for all $i\geq 0$, and $\L_{\geq 1}$ is a finitely
generated left $\L$-module.  Furthermore, assume that all simple
(graded) $\L$-modules or $\L^\op$-module of maximal projective
dimension $n$ are in $\calH^n_\L$ or $\calH^n_{\L^\op}$, respectively,
and if $S$ is a simple (graded) $\L^\op$-module in $\calH^n_{\L^\op}$,
then $\tr^n_{\L^\op}(S)$ has a simple module from $\calH^n_\L$ as a
submodule. Then the category $\calH^n_\L$ consists exactly of the
modules of finite length with composition factors only in $\S^n$.
\end{prop}
\begin{proof}
By Proposition \ref{prop:duality} (a) it follows that all modules of
finite length with composition factors only in $\S^n$ are contained in
$\calH^n_\L$.

Conversely, let $M$ be in $\calH^n_\L$. By Lemma \ref{lem:torsioninH}
we have that $\Ext^i_\L(t^n(M),\L)=(0)$ for $i\neq n$. The exact
sequence $0\to t^n(M)\to M\to M/t^n(M)\to 0$ implies that
$\Ext^i_\L(M/t^n(M),\L)=(0)$ for $i\neq n$. Using Proposition
\ref{prop:torsionfreepddim} it follows that $M/t^n(M)=(0)$. Therefore
$M=t^n(M)$, and since $M$ is finitely generated, $M$ has finite length
and with composition factors only in $\S^n$.
\end{proof}

The above result and the corollary below motivates the following
generalization of Artin-Schelter regular algebras.

\begin{defin}
Let $\L$ be a semiperfect two-sided Noetherian ring, or let $\L$ is a
positively graded ring such that $\L_0\simeq K^t$ for some $t$,
$\dim_K\L_i<\infty$ for all $i\geq 0$, and $\L_{\geq 1}$ is finitely
generated as a left and a right $\L$-module. Assume that $\L$ has
finite (graded) global dimension $n$.

Then $\L$ is \emph{(graded) generalized Artin-Schelter regular of
dimension $n$} if $\calH^n_\L$ and $\calH^n_{\L^\op}$ contain all
simple (graded) $\L$-modules and $\L^\op$-modules of maximal
projective dimension $n$, respectively, and $\tr^n$ maps these simple
modules to modules having a simple module from $\calH^n$ on the
opposite side as a submodule for all (graded) simple $\L^\op$-modules
$S$ of maximal projective dimension $n$.
\end{defin}

We summarize some basic immediate consequences for generalized
Artin-Schelter algebras below.
\begin{cor}
Suppose that $\L$ is (graded) generalized Artin-Schelter regular of dimension
$n$. Then the following assertions hold.
\begin{enumerate}
\item[(a)] The categories $\calH^n_\L$ and $\calH^n_{\L^\op}$ consist
  exactly of the modules of finite length with composition factors
  only in simple (graded) modules of maximal projective dimension $n$.
\item[(b)] The functors $\tr^n_\L\colon \calH^n_\L\to
\calH^n_{\L^\op}$ and $\tr^n_{\L^\op}\colon \calH^n_{\L^\op}\to
\calH^n_\L$ are inverse exact dualities, which preserve length, and
in particular gives rise to a bijection between the simple (graded)
$\L$-modules and simple (graded) $\L^\op$-modules with maximal
projective dimension $n$.
\end{enumerate}
\end{cor}
\begin{proof}
The claim in (a) follows directly from Proposition
\ref{prop:exactlytorsion}, and (b) follows from (a) and Proposition
\ref{prop:duality} (b).
\end{proof}
For a generalized Artin-Schelter regular algebra of dimension $n$ we can
reformulate the definition of the subfunctor $t^n$ of the identity as
follows. Let $M$ be a $\L$-module. Then $t^n(M)=\sum_{L\subseteq M}
L$, where the sum is taken over all finite length submodules $L$ of
$M$ with $L$ in $\calH^n_\L$.

Classically a graded connected $K$-algebra $\L$ is Artin-Schelter
regular if $\gldim\L=n<\infty$, it has finite Gelfand-Kirillov
dimension, and the unique simple graded $\L$-module $S$ satisfies the
Gorenstein condition, that is, $\Ext^n_\L(S,\L)=(0)$ for $i\neq n$ and
$\Ext^n_\L(S,\L)$ is isomorphic to some graded shift of the unique
simple graded $\L^\op$-module. It is then immediate that such a
algebra is also Artin-Schelter regular in our sense.

\section{Properties of generalized Artin-Schelter regular
  algebras}\label{section1.5}

For classical Artin-Schelter regular algebras one has Serre duality,
and there is an intimate relationship with finite dimensional
Frobenius algebras. This section investigates analogues of these
connections for our generalized Artin-Schelter regular algebras, and
we show similar behavior in our case. In particular, if $\L$ is a
(graded) generalized Artin-Schelter regular algebra of dimension $n$
with duality $D$, then
\[D(\Ext^i_\L(M,-))\simeq \Ext^{n-i}_\L(-,D\tr^n_\L M))\]
for all $M$ in $\calH^n_\L$ and all $i$, and the $\Ext$-algebra of any
finite set of simple $\L$-modules i $\calH^n_\L$ permuted by
$D\tr^n_\L$ is Frobenius. We show a partial converse of this in
Section \ref{section:6}.

We begin by discussing Serre duality type formulas for our generalized
Artin-Schelter regular algebras when they admit an ordinary duality.

For generalized Artin-Schelter regular rings $\L$ of dimension $n$,
the category $\calH^n_\L$ is exactly the full subcategory of
$\L$-modules of finite length with composition factors only in
$\S^n$. Hence, $\calH^n_\L$ and $\calH^n_{\L^\op}$ are abelian
categories and $\tr^n_\L\colon \calH^n_\L\to \calH^n_{\L^\op}$ is an
exact duality. If $\L$ is semiperfect two-sided Noetherian and a
$K$-algebra for a field $K$ and all simple $\L$-modules are finite
dimensional, then the duality $D=\Hom_K(-,K)$ induces a duality
between the finite length $\L$-modules and $\L^\op$-modules.  When
$\L$ is graded as usual with $\L_0\simeq K^t$ for some $t$, then the
duality $D=\Hom_K(-,K)$ provides such a duality $D$ as above. When
either of the two cases for a generalized Artin-Schelter regular
algebra $\L$ occurs, $\L$ is said \emph{to have a duality $D$}.
However it is not true that $D$ necessarily induces a duality from
$\calH^n_\L$ to $\calH^n_{\L^\op}$, unless the categories $\calH^n$
consist for all modules of finite length. In this case we have
autoequivalences $D\tr^n_\L$ and $\tr^n_{\L^\op}D$ of $\calH^n_\L$. So
when giving Serre duality type formulas in the next result, the module
$D\tr^n_\L(M)$ is not necessarily in $\calH^n_\L$ again for $M$ in
$\calH^n_\L$.

\begin{prop}\label{prop:formulas}
Let $\L$ be a (graded) generalized Artin-Schelter regular ring of dimension $n$
having a duality $D$. For all pairs of $\L$-modules $L$ and $M$ with
$M$ in $\calH^n_\L$, there are natural isomorphisms
\[\varphi_i\colon D(\Ext^i_\L(M,L))\to \Ext^{n-i}_\L(L,D\tr^n_\L(M))\]
for all $i$, and the Auslander-Reiten formulas
\[D(\Hom_\L(M,L))\simeq \Ext^n_\L(L,D\tr^n_\L(M))\]
and
\[D(\Hom_\L(M,L))\simeq \Ext^1_\L(L,D\Tr(M)).\]
\end{prop}
\begin{proof}
In either situations we have that
$D(\Tor^\L_i(A,B)) \simeq \Ext^i_\L(B,D(A))$ for a $\L^\op$-module $A$,
for a $\L$-module $B$ for all $i\geq 0$. The formula
\[\Tor^\L_i(\tr^n_\L(M),L)\simeq \Ext^{n-i}_\L(M,L)\]
for $M$ in $\calH^n_\L$ from Proposition \ref{prop:duality} induces the
first formula
\[D(\Ext^i_\L(M,L))\simeq \Ext^{n-i}_\L(L,D\tr^n_\L(M))\]
for all $i\geq 0$.  The above formula for $i=0$ gives the first
Auslander-Reiten formula, and using dimension shift the second formula
follows directly.
\end{proof}
Next we show that the $\Ext$-algebra of any finite set of simple
$\L$-modules i $\calH^n_\L$ permuted by $D\tr^n_\L$ for a generalized
Artin-Schelter regular algebra is Frobenius.

Let $\L$ be a generalized Artin-Schelter regular ring of dimension $n$ with a
duality $D$. In our graded situation $\L/\rrad$ with $\rrad$ being the
graded Jacobson radical of $\L$, generates the additive closure of all
simple graded $\L$-modules up to shift. While in the ungraded setting
there might be infinitely many simple modules.  In both cases denote
by $T$ the direct sum of a finite subset of all simple $\L$-modules in
$\S^n$. Denote by $\G$ the Yoneda algebra $\amalg_{i\geq
0}\Ext^i_\L(T,T)$, and consider the natural functors
\[F=\amalg_{i\geq 0}\Ext^i_\L(-,T)\colon \Mod\L \to \Gr(\G)\]
and
\[F'=\amalg_{i\geq 0}\Ext^i_\L(T,-)\colon \Mod\L \to \Gr(\G^\op).\]
Then we can consider the following diagram of functors
\[\xymatrix@C=60pt{
\Mod\L \ar[r]^{F}\ar[d]_{D\tr^n_\L} &
\Gr(\G)\ar[d]^D\\
\Mod\L\ar[r]^{F'} & \Gr(\G^\op)
}\]
The following result shows the above mentioned intimate connection
with Frobenius algebras.
\begin{prop}\label{prop:commdiag}
Let $\L$ be a generalized Artin-Schelter regular ring of dimension $n$ with a
duality $D$. Let $T$ be the direct sum a finite subset of all simple
$\L$-modules in $\S^n$ such that $D\tr^n_\L$ maps $\add T$ to $\add
T$. Denote by $\F(\add T)$ all modules of finite length with
composition factors only from $\add T$.
\begin{enumerate}
\item[(a)] There is a commutative diagram of functors
\[\xymatrix@C=70pt{
\F(\add T) \ar[r]^-{F=\Ext^*_\L(-,T)}\ar[d]_{D\tr^n_\L} &
\gr(\G)\ar[d]^D\\
\F(\add T) \ar[r]^-{F'=\Ext^*_\L(T,-)} & \gr(\G^\op)
}\]
\item[(b)] Let $T$ be such that $D\tr^n_\L$ induces a permutation of
the simples in $\add T$. Then the Yoneda algebra $\G=\amalg_{i\geq
0}\Ext^i_\L(T,T)$ is Frobenius.
\item[(c)] Let $T$ be such that $D\tr^n_\L$ induces a permutation of
the simples in $\add T$. Then there is a commutative diagram of
functors
\[\xymatrix@C=70pt{
\F(\add T) \ar[r]^-{F}\ar[d]_{D\tr^n_\L} &
\gr(\G)\ar[d]^{\Hom_\G(-,\G)D}\\
\F(\add T) \ar[r]^-{F} & \gr(\G)
}\]
\end{enumerate}
\end{prop}
\begin{proof}
(a) We have that $\G_0=\End_\L(T)$ is Noetherian and $\G_i$ is a
finitely generated $\G_0$-module for all $i$, so that $\G$ is a
Noetherian ring (finitely generated $\G_0$-module).  Since the functor
$F$ takes semisimple modules in $\F(\add T)$ to projective modules
(finitely generated modules), $F$ is a half exact functor and $\G$ is
Noetherian, $F|_{\F(\add T)}$ has its image in $\gr(\G)$. Similar
arguments are used for $F'$. Hence the functors in the diagram end up
in the categories indicated.

By Proposition \ref{prop:formulas} the image of the functors $DF$ and
$F'D\tr^n_\L$ on a module $M$ in $\F(\add T)$  are isomorphic as
abelian groups. So we need to show that they in fact are isomorphic as
$\G^\op$-modules. Recall that the isomorphism from $DF(M)$ to
$(F'D\tr^n_\L)(M)$ as abelian groups is induced by the isomorphisms
\[\varphi_i\colon  D(\Ext^{n-i}_\L(M,L))\to \Ext^i_\L(L,D\tr^n_\L(M))\]
for all $i\geq 0$. For any element $\theta$ in $\Ext^j_\L(T,T)$
we want to prove that $\varphi_{i+j}(f)\theta=\varphi_i(f\theta)$ for
all $f$ in $D\Ext^{i+j}_\L(M,T)$.

Represent $\theta$ in $\Ext^j_\L(T,T)$ as a homomorphism
$\theta\colon \Omega^j_\L(T)\to T$. This gives rise to the
following diagram
\begin{small}
\[\xymatrix{
D\Ext^{i+j}_\L(M,T) \ar[r]^{D\Ext^{i+j}_\L(M,\theta)}
\ar[d]^{\varphi_{i+j}} &
D\Ext^{i+j}_\L(M,\Omega^j_\L(T)) \ar[r] \ar[d]^{\varphi_{i+j}} &
D\Ext^i_\L(M,T)\ar[d]^{\varphi_i}\\
\Ext^{n-i-j}_\L(T,D\tr^n_\L(M))
\ar[r]^{\Ext^{n-i-j}_\L(\theta,D\tr^n_\L(M))} &
\Ext^{n-i-j}_\L(\Omega^j_\L(T),D\tr^n_\L(M)) \ar[r] &
\Ext^{n-i}_\L(T,D\tr^n_\L(M))} \]
\end{small}
The first square commutes, since $\varphi_{l}$ is a natural
transformation for all $l$. The horizontal morphisms in the second
square are compositions of connecting homomorphisms or dual
thereof. The isomorphisms $\varphi_l$ are all induced by the
isomorphism for a finitely generated projective module $P$ and any
module $X$
\[\psi\colon D\Hom_\L(P,X)\to \Hom_\L(X,D(P^*))\]
given by $\psi(g)(x)(p^*)=g(p^*(-)x)$ for $x$ in $X$ and $p^*$ in
$P^*$, where $p^*(-)x$ is in $\Hom_\L(P,X)$. This is a natural
isomorphism in both variables. It follows from this that the
isomorphisms $\varphi_l$ all commute with connecting morphisms, that
is, if $0\to X\to Y\to Z\to 0$ is exact, then the diagram
\[\xymatrix{
D\Ext^{i+1}_\L(M,X)\ar[r]^{D(\partial)} \ar[d]^{\varphi_{i+1}} &
D\Ext^i_\L(M,Z)\ar[d]^{\varphi_i}\\
\Ext^{n-i-1}_\L(X,D\tr^n_\L(M))\ar[r]^{\partial'} &
\Ext^{n-i}_\L(Z,D\tr^n_\L(M))}\]
is commutative.

Then going back to our first diagram, which by the above is
commutative, the upper diagonal path is computing $\varphi_i(f\theta)$
while the lower diagonal path is computing $\varphi_{i+j}(f)\theta$
both for any $f$ in $D\Ext^{i+j}_\L(M,T)$. Hence we conclude that
they are equal and $DF(M)$ and $F'(D\tr^n_\L(M))$ are isomorphic as
$\G^\op$-modules for all $M$ in $\F(\add T)$. Since the isomorphism
between $DF$ and $F'D\tr^n_\L$ is induced by the natural isomorphisms
$\{\varphi_i\}_{i=0}^n$, we infer that $DF$ and $F'D\tr^n_\L$ are
isomorphic as functors.

(b) The functor $F\colon \F(\add T)\to \gr(\G)$ sends simple modules
to projective modules, and the functor $D\tr^n_\L\colon \F(\add T)\to
\F(\add T)$ sends simple modules to simple modules. From the
commutative diagram in (a) we infer that the projective module
$\Ext^*_\L(T,D\tr^n_\L(S))$ is isomorphic to the injective module
$D(\Hom_\G(\Ext^*_\L(S,T),\G))$. Since $D\tr^n_\L$ is a permutation of
all simple modules in $\add T$, it follows that the projective and the
injective modules coincide and $\G$ is Frobenius.

(c) Using (a) it is enough to show that $\Hom_\G(-,\G)F'$ is
isomorphic to $F$.

We want to consider $F'$ as a functor $F'\colon E(\F(\add T))\to
\gr(\G)$. Here $E(\F(\add T))$ is the $\Ext$-category of the full
subcategory $\F(\add T)$ of $\mod\L$, that is, $E(\F(\add T))$ has
$\F(\add T)$ as objects and morphisms are given by 
\[\Hom_{E(\F(\add T))}(M,N)=\amalg_{i\geq 0}\Ext^i_\L(M,N)\]
for $M$ and $N$ in $\F(\add T)$. On objects $F'$ is given as usual,
while on morphisms we have the following. For $M$ and $N$ in
$E(\F(\add T))$, we let for $\theta$ in $\Ext^i_\L(M,N)$
\[F'(\theta)\colon F'(M)=\Ext^*_\L(T,M)\to \Ext^*_\L(T,N)=F'(N)\]
be given by the Yoneda multiplication by $\theta$. Hence in this way
$F'$ induces a morphism from $\Ext^*_\L(M,N)$ to
$\Hom_\G(F'(M),F'(N))$. In particular for $N=T$ we have that $F'$
induces a morphism from $F(M)=\Ext^*_\L(M,T)$ to
$\Hom_\G(F'(M),F'(T))=\Hom_\G(F'(M),\G)$.

The above morphism induced by $F'$ is clearly an isomorphism for $T$,
and therefore we have an isomorphism for all simple modules in
$\S^n$. Suppose we have an isomorphism for all modules of length $l$
in $\F(\add T)$. Let $L$ be a module with length $l+1$ in
$\F(\add T)$. Then there is an exact sequence $\eta\colon 0\to L'\to
L\to S\to 0$ with $S$ a simple module in $\S^n$ and $L'$ of length $l$
in $\calH^n_\L$. The long exact sequence induced by $\eta$ applying
$\Hom_\L(T,-)$ and the fact that $\G$ is selfinjective by (b) induce
the following commutative diagram
\[\xymatrix{
(F(L')[-1]\ar[r]\ar[d] &
F(S)\ar[r]\ar[d]     &
 F(L)\ar[r]\ar[d]    &
F(L')\ar[r]\ar[d]    &
F(S)[1]\ar[d]      \\
 F'(L')[1])^*\ar[r] &
(F'(S))^*\ar[r] &
(F'(L))^*\ar[r] &
(F'(L'))^*\ar[r] &(F'(S)[-1])^*  }\]
Since by induction the two leftmost and the two rightmost vertical
maps are isomorphisms, we have by the Five Lemma that $F(L)\simeq
\Hom_\G(\Ext^*_\L(T,L),\G)$. The claim follows from this.
\end{proof}

For connected Koszul algebras $\L$ of finite global dimension P.\
Smith in \cite{S} has shown that the Koszul dual of $\L$ is Frobenius
if and only if $\L$ is Gorenstein. This result is generalized to the
non-connected situation by Martin\'ez-Villa in \cite{MV}.
Lu-Palemieri-Wu-Zhang  have shown in \cite{LPWZ} that for a connected
graded algebra $\L$ with $\Ext$-algebra $E(\L)$, then $\L$ is
Artin-Schelter regular if and only if $E(\L)$ is Frobenius. The above
result shows one direction of this result for our more general class
of Artin-Schelter regular algebras. For illustration let us return to
Example \ref{exam:Auslanderalgebra}.

\begin{example}\label{exam:AuslanderalgebraII}
Let $\Sigma=K[x]/(x^3)$ as in Example \ref{exam:Auslanderalgebra}, and
let $\L$ be the Auslander algebra $\End_\Sigma(M)^\op$ of $\Sigma$
where $M=U\amalg L\amalg \Sigma$. Applying the functor
\[\Hom_\L(-,\L)=\Hom_\L(-,\Hom_\Sigma(M,M)^\op)\] to the
projective resolutions of $S_U$ and $S_L$ induces the following exact
sequences
\[0\to {_\Sigma(U,M)}\to {_\Sigma(L,M)}\to {_\Sigma(U,M)}\to
\tr^2_\L(S_U)\to 0\]
and
\[0\to {_\Sigma(L,M)}\to {_\Sigma(\Sigma\amalg U,M)}\to {_\Sigma(L,M)}\to
\tr^2_\L(S_L)\to 0,\]
where $\tr^2_\L(S_U)\simeq S_U^\op$ and $\tr^2_\L(S_L)\simeq
S_L^\op$. Since $D(S^\op)\simeq S$ for all simple module $S$ (and
$S^\op$ over $\L^\op$), it follows that $D\tr^2_\L(S)\simeq S$ for
$S=S_U, S_L$. Hence we see that the condition in Proposition
\ref{prop:commdiag} (b) is satisfied, and therefore $\G=\Ext^*_\L(T,T)$ is
Frobenius for $T=S_U\amalg S_L$, as pointed out already in Example
\ref{exam:Auslanderalgebra}.

For an Auslander algebra $\L$ of a finite dimensional algebra $\Sigma$
of finite representation type, it is straightforward to see that
$D\tr^2_\L(S_C)$ for some indecomposable $\Sigma$-module $C$ is given
by $S_{\tau C}$, where $\tau$ is the Auslander-Reiten translate. Hence
to find a finite set of simple modules satisfying the assumptions in
Proposition \ref{prop:commdiag} (b), it is enough to choose simple
modules corresponding to indecomposable $\Sigma$-modules permuted by
$\tau$.
\end{example}

Now we make a further remark concerning Proposition
\ref{prop:commdiag}. Recall that the functor $\Hom_\G(-,\G)D$ is the
Nakayama functor for the Frobenius algebra $\G$. Hence, Proposition
\ref{prop:commdiag} (c) in other words says that $D\tr^n_\L$
essentially behaves like the Nakayama functor via the functor
$F=\Ext^*_\L(-,T)$.

For hereditary categories discussed in \cite{RvdB}, the existence of a
Serre functor is connected with the existence of almost split
sequences. In Proposition \ref{prop:formulas} we saw that $D\tr^n_\L$
plays the role of a Serre functor for generalized Artin-Schelter
regular algebras. For our $n$-dimensional Artin-Schelter regular
algebras, it gives rise to $n$-fold almost split extensions, as we
explain next. To this end recall that a non-zero $n$-fold extension
$\eta$ in $\Ext^n_\L(C,A)$ is called an \emph{$n$-fold almost split
  extension} if for any non-splittable epimorphism $f\colon X\to C$
and for any non-splittable monomorphism $g\colon A\to Y$ the pullback
$\eta\cdot f$ and the pushout $g\cdot \eta$ are both zero in
$\Ext^n_\L(X,A)$ and $\Ext^n_\L(C,Y)$, respectively. With these
remarks we give a connection between the Serre functor $D\tr^n_\L$ and
existence of $n$-fold almost split extensions for the category
$\calH^n_\L$ for generalized Artin-Schelter regular algebras.

\begin{prop}\label{prop:almostsplit}
Let $\L$ be a generalized Artin-Schelter regular ring of dimension $n$
with duality $D$. Assume that $\calH^n_\L$ is a Krull-Schmidt category.
\begin{enumerate}
\item[(a)] If $n=1$, then the category $\calH^1_\L$ has almost split
  sequences.
\item[(b)] If $n>1$, then the category $\calH^n_\L$ has $n$-fold
  almost split extensions.
\end{enumerate}
\end{prop}
\begin{proof}
From Proposition \ref{prop:formulas} we have that
\[D(\End_\L(M))\simeq \Ext^n_\L(M,D\tr^n_\L(M))\]
for any module $M$ in $\calH^n_\L$. If $M$ is indecomposable, then
$\End_\L(M)$ is a local ring. Then using standard arguments as in
\cite{AR} both claims in (a) and (b) follows.
\end{proof}

\section{The polynomial ring}\label{section2}

In this section we consider the polynomial ring $\L=K[x_1,\ldots,x_n]$
in $n$ indeterminants $x_1,\ldots,x_n$ over an algebraically closed
field $K$. The aim is to illustrate some of the aspects discussed in
the previous section on this concrete example.

First we show that $\L$ is a generalized Artin-Schelter regular ring
of dimension $n$ as a ungraded ring.  The simple modules over $\L$ are
well-known.  Since $K$ is algebraically closed, all maximal ideals are
of the form
\[\frakm_{\overline{a}}=(x_1-a_1,x_2-a_2,\ldots,x_n-a_n)\]
for some $n$-tuple $\overline{a}=(a_1,a_2,\ldots, a_n)$ in
$K^n$. Hence all simple $\L$-modules are given as
$S_{\overline{a}}=\L/\frakm_{\overline{a}}$.

For any $n$-tuple $\overline{a}=(a_1,a_2,\ldots, a_n)$ in $K^n$ there
is an automorphism $\sigma_{\overline{a}}\colon \L\to \L$ given by
$\sigma_{\overline{a}}(x_i)=x_i-a_i$ for all $i=1,2,\ldots,n$.

For a $\L$-module $M$ and for an automorphism $\sigma\colon \L\to \L$
denote by $M^\sigma$ the $\L$-module with underlying vector space $M$
and where the action of an element $\lambda$ in $\L$ is given by
$\lambda\cdot m = \sigma^{-1}(\lambda)m$. This extends to a functor
$\sigma\colon \Mod\L\to \Mod\L$, which is an exact functor. Since the
functor $\Hom_\L(X^\sigma,-)\simeq \Hom_\L(X,-)\sigma^{-1}(-)$ for an
automorphism of $\L$, the functor $\sigma$ preserves projective
modules and in particular free modules.

With the above convention it follows that $S_{\overline{a}}\simeq
K^{\sigma_{\overline{a}}}$ for all $n$-tuples $\overline{a}$ in
$K^n$. In particular we have that $K\simeq S_{0}$ with
$0=(0,0,\ldots,0)$.

Let
\[\mathbb{P}\colon 0\to P_n\to P_{n-1}\to \cdots \to P_1\to P_0\to
S_0\to 0\] be a minimal graded projective resolution (the Koszul
complex) of $S_0=K$ over $\L$. This is a projective resolution of
$S_0$ also when we forget the grading, and we furthermore have that
$\Ext^i_\L(S_0,\L)=(0)$ for $i\neq n$ and $\Ext^n_\L(S_0,\L)\simeq
S_0$. Let $\overline{a}=(a_1,a_2,\ldots,a_n)$ be in $K^n$. Then
$\mathbb{P}^{\sigma_{\overline{a}}}$ is a projective resolution of
$S_{\overline{a}}$. Since $\L^{\sigma_{\overline{a}}}\simeq \L$ for
all $n$-tuples $\overline{a}$ and $\Hom_\L(X^\sigma,\L)\simeq
\Hom_\L(X,\L^{\sigma^{-1}})$, we infer that
$\Ext^i_\L(S_{\overline{a}}, \L)=(0)$ for all $i\neq n$ and
$\tr^n_\L(S_{\overline{a}})=\Ext^n_\L(S_{\overline{a}},\L)\simeq
S_{\overline{a}}$.  Consider as before the subcategory
\[\calH^n_\L=\{ M\in\mod\L\mid \pd_\L M=n, \Ext^i_\L(M,\L)=(0) \text{\
  for all\ } i\neq 0\}.\]
By the above observations all simple $\L$-modules are in $\calH^n_\L$
and $\tr^n_\L$ is the identity ($\L=\L^\op$), hence $\L$ is a
generalized Artin-Schelter regular $K$-algebra as a ungraded ring.  As
a consequence of the above and Proposition \ref{prop:exactlytorsion}
we have the following.
\begin{thm}
Let $\L=K[x_1,x_2,\ldots,x_n]$ be the polynomial ring in $n$
indeterminants $x_1,x_2,\ldots,x_n$ over an algebraically closed field
$K$. Then $\calH^n_\L$ is exactly the full subcategory consisting of
all $\L$-modules of finite length in $\Mod\L$.
\end{thm}

Now let us discuss Proposition \ref{prop:commdiag} (b) for this
example. In particular it says that the Yoneda algebra $\amalg_{i\geq
  0}\Ext^i_\L(K,K)$ is Frobenius. This is well-known and it is in fact
isomorphic to the exterior algebra on $K^n$.  However we can consider
$\L$ as an ungraded $K$-algebra, and we have seen that $\L$ is a
generalized Artin-Schelter regular $K$-algebra. Consider two simple
$\L$-modules $S_{\overline{a}}$ and $S_{\overline{b}}$. Then each
$\Ext^i_\L(S_{\overline{a}}, S_{\overline{b}})$ is a module over $\L$,
where the action of $\L$ can be taken through $S_{\overline{a}}$ or
$S_{\overline{b}}$. Hence
$\Ext^i_\L(S_{\overline{a}},S_{\overline{b}})$ is annihilated by
$\frakm_{\overline{a}}+\frakm_{\overline{b}}=\L$ if $\overline{a}\neq
\overline{b}$.  Therefore the $\Ext$-algebra of any finite number of
simple modules over $\L$ is just a direct sum of the $\Ext$-algebras
of the simple $\L$-modules involved.  Since $S_{\overline{a}} \simeq
K^{\sigma_{\overline{a}}}$, it follows that each of these
$\Ext$-algebras are isomorphic to the exterior algebra. And hence
Frobenius as predicted by Proposition \ref{prop:commdiag} (b).

Finally we point out that $\calH^n_\L$ as $n$-fold almost split
sequences.  Since $\Ext^i_\L(S_{\overline{a}},S_{\overline{b}})=(0)$
for all $i\geq 0$ whenever $\overline{a}\neq \overline{b}$, any object
of finite length in $\calH^n_\L$ is a module over
$\L/{\mathfrak{m}_{\overline{a}}}^n$ for some $\overline{a}$ and some
positive integer $n$. Since this is a finite dimensional algebra, it
follows that $\calH^n_\L$ is a Krull-Schmidt category. As a corollary
of Proposition \ref{prop:almostsplit} we have the following.

\begin{cor}
Let $\L=K[x_1,x_2,\ldots,x_n]$ be the polynomial ring in $n$
indeterminants over an algebraically closed field $K$. Then the
category of $\L$-modules of finite length $\calH^n_\L$ has $n$-fold
almost split sequences.
\end{cor}

\section{Non-Noetherian generalized Artin-Schelter regular
  algebras}\label{section3}\

It is unknown whether or not (to our knowledge) that all connected
Artin-Schelter regular algebras are Noetherian.  In our definition a
generalized Artin-Schelter regular algebra we have not put any
requirement on the Gelfand-Kirillov dimension of of the algebra. So
this is maybe why we in this section can exhibit examples of
generalized Artin-Schelter regular algebras which are not
Noetherian. However, they are Noetherian if and only if the
Gelfand-Kirillov dimension is finite.

To give the class of generalized Artin-Schelter regular algebras we
construct in this section, we recall for the convenience of the reader
the results we used from \cite{BGL,GMT,M1,M2,M3}. The first result we
recall characterizes selfinjective Koszul algebras in terms of the
Koszul dual. 
\begin{thm}[\cite{M1}]
  Let $\Lambda$ be a Koszul algebra with Yoneda algebra $\Gamma$. Then
  $\Lambda$ is selfinjective if and only if there exists some positive
  integer $n$ such that all graded simple $\Gamma$-modules have the
  same projective dimension $n$ and $\Gamma$ satisfies the $n$-simple
  condition.
\end{thm}
It is easy to find selfinjective Koszul algebras, as the following theorem
shows.
\begin{thm}[\cite{M2}]
  Let $\Lambda$ be an indecomposable selfinjective finite dimensional
  $K$-algebra with radical $\rrad$ such that $\rrad^{2}\neq 0$ and
  $\rrad^{3}=0$. Then $\Lambda$ is Koszul if and only if it is of
  infinite representation type. In the Koszul case the Yoneda algebra
  $\Gamma$ has global dimension $2$.
\end{thm}

\begin{cor}
  Let $Q$ be a connected bipartite graph, and let $K$ be a field. Then
  the trivial extension $\Lambda =KQ\ltimes D(KQ)$ is Koszul
  if and only if $Q$ is non-Dynkin. If $Q$ is non-Dynkin, then the
  Yoneda algebra $\Gamma$ of $\Lambda$ is the preprojective algebra.
\end{cor}

We can characterize when an Artin-Schelter regular Koszul algebra of
global dimension $2$ is Noetherian.

\begin{thm}[\cite{BGL,GMT}]\label{thm:bgl+gmt}
  Let $\Gamma$ be an indecomposable Koszul algebra, where all the
  graded simple $\Gamma$-modules have projective dimension $2$ and
  satisfying the $2$-simple condition. Suppose $\L$ is the Yoneda
  algebra of $\Gamma$. If $\Lambda$ is a tame algebra, then $\Gamma$
  is Noetherian of Gelfand-Kirillov dimension two. If $\Lambda$ is
  wild, then $\Gamma$ is non-Noetherian of infinite Gelfand-Kirillov
  dimension.
\end{thm}

It is easy to built new Artin-Schelter regular algebras from given ones,
this is the situation considered in next theorem.

\begin{thm}[\cite{M3}]\label{thm:m3}
  Let $\Lambda$ be a Koszul $K$-algebra with Yoneda algebra $\Gamma$,
  and let $G$ be a finite group of automorphisms of $\Lambda$, such
  that characteristic of $K$ does not divide the order of $G$.  Then
  the following statements are true.
\begin{enumerate}
\item[(a)] There is a natural action of $G$ on $\Gamma$.
\item[(b)] The skew group algebra $\Lambda \ast G$ is Koszul with
  Yoneda algebra $\Gamma \ast G$.
\item[(c)] If $\Lambda$ is selfinjective, then $\Lambda \ast G$ is
  selfinjective, in particular all the graded simple $\Gamma\ast
  G$-modules have projective dimension $n$ and $\Gamma \ast G$
  satisfies the $n$-simple condition, where $n$ is the Loewy length of
  $\L\ast G$. 
\end{enumerate}
\end{thm}
Recall that when $\L$ is a graded algebra and $G$ is a group of
automorphisms acting on $\L$, then $\L\ast G$ is a graded algebra
again with degree $i$ of $\L\ast G$ being given by $\L_i\ast G$.

\begin{cor}
  Let $\Gamma = K[ x_{1},x_{2},\dots,x_{n}] $ be the polynomial algebra,
  and let $G$ be a finite group of automorphisms of $\Gamma $ such
  that characteristic of $K$ does not divide the order of $G$. Then
  the skew group algebra $K[ x_{1},x_{2},\ldots,x_{n}] \ast G$ is
  Artin-Schelter regular.
\end{cor}
\begin{proof}
  Let $\Lambda$ be the exterior algebra of a $K$-vector space of
  dimension $n$. Then it is known $\Lambda$ is selfinjective Koszul
  with Yoneda algebra $\Gamma$. Then apply Theorem \ref{thm:m3}, and
  the claim follows.
\end{proof}

Now we are ready to give the construction of a class of generalized
Artin-Schelter regular algebras in the sense of Section
\ref{section:1}. Let $Q$ be the quiver
$\xymatrix{1\ar@<2ex>[r]^{\alpha_1}\ar@<1ex>[r]_{\alpha_2}
  \ar@{}@<2.6ex>[r]_(.3){\vdots}\ar@<-2ex>[r]_{\alpha_n} & 2}$ for
$n\geq 2$, where we denote the trivial paths corresponding to the
vertices by $e_1$ and $e_2$. Consider the trivial extension $\Lambda =
KQ\vartriangleright D(KQ)$. Then $\Lambda$ is isomorphic to the
algebra $K\widehat{Q}/I$, where $ \widehat{Q}$ is the quiver described
as follows. The quiver $\widehat{Q}$ has vertices
$\widehat{Q}_{0}=Q_{0}$ and arrows $\widehat{Q}_{1}=Q_{1}\cup
Q_{1}^\op$, that is, in addition to the old arrows in $Q$, we
introduce an arrow $\widehat{\alpha}$ in the opposite direction for
any arrow $\alpha$ in $Q$. The ideal $I$ is generated by the following
relations $\{\alpha\widehat{\alpha}-\beta\widehat{ \beta},
\widehat{\alpha}\alpha-\widehat{\beta}\beta, \widehat{\alpha}\beta,
\alpha\widehat{\beta}\}_{\alpha\neq \beta\in Q_1}$.

Let the cyclic group $G=\langle g\rangle$ of order $2$ act as automorphisms
of $\Lambda $ by
\begin{align}
g(e_{1}) & = e_{2},\notag\\
g(e_{2}) & = e_{1},\notag\\
g(\alpha) & = \widehat{\alpha},\notag\\
g(\widehat{\alpha}) & = \alpha.\notag
\end{align}
Assume the characteristic of $K$ different from 2. Then $\Lambda \ast
G$ is selfinjective Koszul, and $(\Lambda \ast G)_{0}$ has
$\{\overline{e}_{i}=e_{i}\otimes 1\}_{i=1}^2\cup
\{\overline{e}_{i}g=e_{i}\otimes g\}_{i=1}^2$ as a basis over
$K$. Define the map $\varphi\colon (\L\ast G)_0\to M_2(K)$ by letting
\begin{align}
\varphi(\overline{e}_{1}) & = \left(
\begin{smallmatrix}
1 & 0 \\
0 & 0%
\end{smallmatrix}%
\right),\notag\\
\varphi(\overline{e}_{2}) & = \left(
\begin{smallmatrix}
0 & 0 \\
0 & 1%
\end{smallmatrix}%
\right) ,\notag\\
\varphi(\overline{e}_{1}g) & = \left(
\begin{smallmatrix}
0 & 1 \\
0 & 0%
\end{smallmatrix}%
\right),\notag \\
\varphi(\overline{e}_{2}g) & = \left(
\begin{smallmatrix}
0 & 0 \\
1 & 0%
\end{smallmatrix}%
\right).\notag
\end{align}
One easily checks that $\varphi$ is an isomorphism and therefore
$(\Lambda \ast G)_{0}\simeq M_2(K)$.

The degree one part $(\Lambda\ast G)_1$ is generated by
$\langle\alpha\otimes 1,\alpha\otimes g,\widehat{\alpha}\otimes 1,
\widehat{\alpha}\otimes g\rangle_{\alpha\in Q_1}$, and $(\Lambda\ast
G)_2$ is generated by $\langle\alpha\widehat{\alpha}\otimes 1,
\alpha\widehat{\alpha}\otimes g, \widehat{\alpha}\alpha\otimes 1,
\widehat{\alpha}\alpha\otimes g\rangle_{\alpha\in Q_1}$.
It follows that $\L\ast G$ is Morita equivalent to $e_1(\L\ast
G)e_1$. The algebra $e_1(\L\ast G)e_1$ is a basic connected
selfinjective Koszul algebra, where
\begin{align}
(e_{1}(\Lambda \ast G)e_{1})_{0} & = K,\notag\\
(e_{1}(\Lambda \ast G)e_{1})_{1} & =
 \langle\widehat{\alpha}\otimes g\rangle_{\alpha\in Q},\notag\\
(e_{1}(\Lambda \ast G)e_{1})_{2} & =\langle\widehat{\alpha}\alpha\otimes
1\rangle_{\alpha\in Q},\notag
\end{align}
since $e_{1}\widehat{\alpha}\otimes
ge_{1}=e_{1}\widehat{\alpha}e_{2}\otimes g$ and
$(\widehat{\alpha}\otimes g)(\widehat{\alpha}\otimes g) =
\widehat{\alpha}\alpha\otimes g^{2} = \widehat{\alpha}\alpha\otimes 1$.
It follows from this that $e_{1}(\Lambda \ast G)e_{1}$ is isomorphic to
\begin{align}
A & = K \langle
x_{1},x_{2},\ldots,x_{n}\rangle/\langle \{x_{i}x_{j}\}_{i\neq
  j},\{x_{i}^{2}-x_{2}^{2}\}_{i\neq 2}\rangle\notag\\
  & = K[ x_{1},x_{2},\ldots,x_{n}]
/(\{x_{i}x_{j}\}_{i\neq j},\{x_{i}^{2}-x_{2}^{2}\}_{i\neq 2}).\notag
\end{align}
The algebra $A$ is selfinjective Koszul with Yoneda algebra
\[B=A^{!}=
K\langle x_{1},x_{2},\ldots,x_{n}\rangle/\langle\sum_{i=1}^n
x_{i}^{2}\rangle.\]
Since $B$ is Morita equivalent to $\Gamma \ast G$, with $\Gamma $ the
preprojective algebra of the quiver $Q\colon
\xymatrix{1\ar@<2ex>[r]^{\alpha_1}\ar@<1ex>[r]_{\alpha_2}
  \ar@{}@<2.6ex>[r]_(.3){\vdots}\ar@<-2ex>[r]_{\alpha_n} & 2,}$ then
it is non-Noetherian and of Gelfand-Kirillov dimension infinite for
$n\geq 3$ by Theorem \ref{thm:m3} and \ref{thm:bgl+gmt}. We obtain the
following theorem.
\begin{thm}
  The algebras $\Gamma_n = K\langle
  x_{1},x_{2},\ldots,x_{n}\rangle/\langle\sum_{i=1}^nx_{i}^{2}\rangle$
  are Koszul for all $n\geq 1$, and all graded simple
  $\Gamma_n$-modules have projective dimension $n$ and $\Gamma_n$
  satisfies the $n$-simple condition.  For $n=2$ they are Noetherian of
  Gelfand-Kirillov dimension $2$, and for $n\geq 3$ they are
  non-Noetherian and of infinite Gelfand-Kirillov dimension.
\end{thm}

\section{Generalized Artin-Schelter regular
  categories}\label{section4}

The aim of this section is to generalize the notion of generalized
Artin-Schelter regular algebras introduced in Section \ref{section:1}
to generalized Artin-Schelter regular categories. Furthermore, we show
that the same results hold true in this situation. The need for such a
generalization was already indicated in the introduction and Section
\ref{section:1}, as our main application is the associated graded
Auslander category of a component in the Auslander-Reiten quiver of a
finite dimensional algebra $\L$ (See the introduction to Section
\ref{section:7} for further details). In this category there are as
many simple objects as there are indecomposable finitely generated
modules over $\L$. Hence the need to extend to categories is apparent.

Now we proceed with describing the categorical version of generalized
Artin-Schelter regular algebra, where we refer the reader to
\cite{MVS1} for the details concerning the notion discussed below. Let
$\C$ be an additive $K$-category over a field $K$. Consider the
category $\Mod(\C)$ of all additive functors $\C^\op\to \Mod K$, where
the finitely generated projective objects are given by $\Hom_\C(-,C)$
for an object $C$ in $\C$. Define the functor $(-)^*\colon \Mod(\C)\to
\Mod(\C^\op)$ as follows. For $F$ in $\Mod(\C)$ let
\[F^*(X)=\Hom_{\Mod(\C)}(F,\Hom_\C(-,X))\]
for all objects $X$ in $\C$, and for a morphism $f\colon X\to Y$ in
$\C$ we have
\[F^*(f)\colon F^*(X)=\Hom_\C(F,\Hom_\C(-,X))\to
\Hom_\C(F,\Hom_\C(-,Y))=F^*(Y)\]
given by $\Hom_\C(F,\Hom_\C(-,f))$. For a morphism $\eta\colon F\to
F'$ in $\Mod(\C)$ let
\[\eta^*(X)=\Hom_{\Mod(\C)}(\eta,\Hom_\C(-,X))\colon (F')^*(X)\to
F^*(X).\]
In particular we see that $\Hom_\C(-,C)^*=\Hom_\C(C,-)$ for all
objects $C$ in $\C$. Furthermore, for any $F$ in $\Mod(\C)$ there is an
isomorphism
\[\Hom_\C(-,C)^*\otimes_\C F\simeq \Hom_{\Mod(\C)}(\Hom_\C(-,C),F)\]
for any object $C$ in $\C$.

Define for any integer $i\geq 0$ a functor $\tr^i_\C\colon \Mod(\C)\to
\Mod(\C^\op)$ as follows. For $F$ in $\Mod(\C)$ let
\[\tr^i_\C(F)(X)=\Ext^i_{\Mod(\C)}(F,\Hom_\C(-,X))\]
for any object $X$ in $\C$. The action on morphisms in $\C^\op$ and in
$\Mod(\C)$ are defined in the natural way.

Let $n$ be a positive integer. Denote by $\calH^n_\C$ the full
subcategory of $\Mod(\C)$ consisting of all objects $F$ with projective
dimension $n$, a projective resolution consisting of finitely
generated functors in $\Mod(\C)$ and $\tr^i_\C(F)=0$ for all $i\neq n$,

We constantly switch between a graded and a ungraded category $\C$,
and to be more compact in the presentation we use the convention that
if $\C$ is a graded category, then a $\C$-module is an object in
$\Gr(\C)$, and if $\C$ is a ungraded category, then a $\C$-module is
an object in $\Mod(\C)$.

With these preliminaries we have the analogue of Proposition
\ref{prop:duality} from the algebra situation. The proof is literary
the same and it is left to the reader.
\begin{prop}\label{prop:duality2}
\begin{enumerate}
\item[(a)] \sloppy The subcategory $\calH_\C^n$ is closed under
extensions. If the category of $\C$-modules has global dimension $n$,
then $\calH_\C^n$ is closed under cokernels of monomorphisms.
\item[(b)] The functor $\tr^n_\C\colon \Mod(\C)\to \Mod(\C^\op)$ restricts
to a functor $\tr^n_\C\colon \calH^n_\C\to \calH^n_{\C^\op}$, and here
it is an exact duality.
\item[(c)] For any object $M$ in $\calH^n_\C$ and any $\C$-module $L$
there are isomorphisms
\[\Tor^{\Mod(\C)}_i(\tr^n_\C(M),L) \simeq \Ext^{n-i}_{\Mod(\C)}(M,L)\]
for all $i\geq 0$.
\end{enumerate}
\end{prop}
Similarly as for algebras, $\tr^n_{\C^\op}(S)$ is not necessarily a
simple object again, but the following lemma give a sufficient
condition to ensure this. The proof is the same as for algebras, and
it is left to the reader.
\begin{lem}\label{lem:simpletosimplecategory}
Let $\C$ be an additive $K$-category with finite global dimension $n$,
and let $S$ be a simple object in $\calH^n_{\C^\op}$. Assume that
$\tr^n_{\C^\op}(S)$ contains a simple object in $\calH^n_\L$ as a
subobject. Then $\tr^n_{\C^\op}(S)$ is a simple object (and it is in
$\calH^n_\C$).
\end{lem}
To proceed we need to put further conditions on the categories. In
particular, information about the simple objects and projective
covers. If $\C$ is a positively graded Krull-Schmidt $K$-category with
$\rad(-,-)=\amalg_{i\geq 1} \Hom_\C(-,-)_i$, we observed in
\cite[Lemma 2.3]{MVS1} that any bounded below functor $F$ has a
projective cover $P\to F$ with $P/\rad P\simeq F/\rad F$. If $\C$ is a
(positively graded) Krull-Schmidt $K$-category, we showed in
\cite[Lemma 2.4]{MVS1} that the all simple objects are exactly those
of the form $\Hom_\C(-,C)/\rad\Hom_\C(-,C)$ for an indecomposable
object $C$ in $\C$.

When $\C$ is a Krull-Schmidt $K$-category, then the full subcategory
$\mod\C$ of $\Mod(\C)$ consisting of the finitely presented objects
all have projective covers. In addition if $P\to F$ is a projective
cover of a finitely presented object $F$, then $P/\rad P\simeq F/\rad
F$ as above. Moreover we have Nakayama's Lemma in this setting.

Let $\S^n$ denote the simple $\C$-modules which are in
$\calH^n_\C$. Then $\calH^n_\C$ contains all objects of finite length
with composition factors only in $\S^n$ by Proposition
\ref{prop:duality2}. We continue to mimic the algebra situation, and
next we want to find sufficient conditions on $\C$ such that
$\calH^n_\C$ is exactly the full category of $\C$-modules with finite
length and composition factors only in $\S^n$.

Define a subfunctor $t^n\colon \Mod(\C)\to \Mod(\C)$ of the identity
functor as follows. Let $M$ be a $\C$-module. Then let
$t^n(M)=\sum_{L\subseteq M} L$, where the sum is taken over all
subobjects $L$ of $M$, where $L$ has finite length with composition
factors only in $\S^n$. We say that a $\C$-module $M$ is
\emph{$\calH^n_\C$-torsion} if $M=t^n(M)$. A $\C$-module $M$ is
\emph{$\calH^n_\C$-torsion free} if $t^n(M)=(0)$.

\begin{lem}\label{lem:torsionradical}
Let $\C$ be a (positively graded) Krull-Schmidt (locally finite)
$K$-category such that all simple $\C$-modules are finitely presented.
Then $t^n(M/t^n(M))=(0)$ for all $\C$-modules $M$, or
equivalently $M/t^n(M)$ is $\calH^n_\C$-torsion free.
\end{lem}
\begin{proof}
Let $M$ be a $\C$-module and assume that $M/t^n(M)$ has a subfunctor
of finite length with composition factors only in $\S^n$. Then
$M/t^n(M)$ has a simple subobject $S$ from $\S^n$. Let $S'$ be the
pullback of $S$ in $M$, so that we have an exact sequence
\[0\to t^n(M)\to S'\to S\to 0.\]
It induces the following commutative exact diagram.
\[\xymatrix{
0\ar[r]  & \rad_\C(-,C)\ar[r]\ar[d]^h
         & (-,C)\ar[r]\ar[d]^f
         & S\ar[r]\ar@{=}[d] & 0\\
0\ar[r]  & t^n(M)\ar[r] & S'\ar[r] & S\ar[r] & 0}\]
Since $\rad_\C(-,C)$ is finitely generated, $\Im h$ is of finite
length with composition factors in $\mathcal{S}^{n}$. We infer that
$\Im f$ has finite length with composition factors only in $\S^n$, so
that $\Im f\subseteq t^n(M)$. This is a contradiction, so that
$t^n(M/t^n(M))=(0)$.
\end{proof}
As for modules, we need that the finitely generated
$\calH^n_\C$-torsion free objects have projective dimension at most
one less than the global dimension, as we show next.
\begin{lem}\label{lem:projdimtorfree}
Let $\C$ be a (positively graded) Krull-Schmidt (locally finite)
$K$-category such that the category of $\C$-modules has finite global
dimension $n$ and that all simple $\C$-modules are finitely
presented. Suppose that $\Mod(\C)$ is a Noetherian category, or when $\C$
is graded that $\rad =\Hom_\C(-,-)_{\geq 1}$. Suppose that all simple
$\C^\op$-modules of maximal projective dimension are in
$\calH^n_{\C^\op}$ and that $\tr^n_{\C^\op}(S)$
contains a simple $\C$-module from $\calH^n_\C$ as a subobject for all
simple $\C^\op$-modules $S$ of maximal projective dimension $n$. Then any
finitely generated $\calH^n_\C$-torsion free $\C$-module has
projective dimension at most $n-1$.
\end{lem}
\begin{proof}
Let $L$ be a finitely generated torsion free $\C$-module. As
for modules we obtain
$\Tor^{\Mod(\C)}_n(\Hom_\C(C,-)/\rad\Hom_\C(C,-),L)=(0)$ for all
indecomposable objects $C$ in $\C$. Hence we infer that
$\Omega^n_{\C}(L)/\rad\Omega^n_\C(L)=(0)$. By Nakayama's Lemma
\cite[Lemma 1.10]{MVS1} we have
that $\Omega^n_\C(L)=(0)$, and the projective dimension of $L$ is at
most $n-1$.
\end{proof}

The proof of the following lemma is also the same as the one for the
module situation. It shows that any $\calH^n_\C$-torsion $\C$-module
$M$ satisfies the same $\Ext$-vanishing conditions as $\C$-modules in
$\calH^n_\C$, namely $\tr^i_\C(M)=(0)$ for all $i\neq n$.
\begin{lem}
Let $\C$ be a (positively graded) Krull-Schmidt (locally finite)
$K$-category such that the category of $\C$-modules has finite global
dimension $n$ and all simple $\C$-modules are finitely
presented. Suppose that all simple $\C$-modules of maximal projective
dimension are in $\calH^n_\C$. Then any $\calH^n_\C$-torsion
$\C$-module $M$ satisfies $\tr^i_\C(M)=(0)$ for all $i\neq n$.
\end{lem}

The fact that $\calH^n_\C$ consists exactly of the objects of finite
length now follows as for modules.
\begin{prop}
Let $\C$ be a (positively graded) Krull-Schmidt (locally finite)
$K$-category such that the category of $\C$-modules has finite global
dimension $n$ and that all simple $\C$-modules are finitely
presented. Suppose that $\Mod(\C)$ is Noetherian category, or when $\C$
is graded $\rad =\Hom_\C(-,-)_{\geq 1}$. Suppose that all simple
$\C$-modules and $\C^\op$-modules of maximal projective dimension are
in $\calH^n_{\C^\op}$ and $\calH^n_{\C^\op}$, respectively, and if $S$
is a simple $\C^\op$-module in $\calH^n_{\C^\op}$, then
$\tr^n_{\C^\op}(S)$ contains a simple $\C$-module from $\calH^n_\C$ as
a subobject for all simple $\C^\op$-modules $S$ of maximal projective
dimension $n$.

Then the category $\calH^n_\L$ consists exactly of the $\C$-modules of
finite length with composition factors only in $\S^n$.
\end{prop}
This leads to the following definition of a generalized Artin-Schelter
regular category.
\begin{defin}
Let $\C$ be a (positively graded) Krull-Schmidt (locally finite)
$K$-category such that the category of $\C$-modules has finite global
dimension $n$. Let $\C$ be a Noetherian category, or let $\C$ be a
positively graded such that $\rad =\Hom_\C(-,-)_{\geq 1}$.

Then $\C$ is \emph{generalized Artin-Schelter regular (of dimension
  $n$)} if $\calH^n_\C$ and $\calH^n_{\C^\op}$ contain all simple
$\C$-modules and $\C^\op$-modules of maximal projective dimension $n$,
respectively, and $\tr^n$ maps these simple $\C$-modules to objects
having a simple $\C$-module from $\calH^n$ on the opposite side as a
subobject.
\end{defin}
The basic elementary properties of generalized Artin-Schelter regular
categories are the following.
\begin{cor}
Suppose that $\C$ is generalized Artin-Schelter regular $K$-category
of dimension $n$. Then the following assertions hold.
\begin{enumerate}
\item[(a)] The categories $\calH^n_\C$ and $\calH^n_{\C^\op}$ consist
exactly of the objects of finite length with composition factors only
in simple $\C$-modules of maximal projective dimension $n$.
\item[(b)] The functors $\tr^n_\C\colon \calH^n_\C\to
\calH^n_{\C^\op}$ and $\tr^n_{\C^\op}\colon \calH^n_{\C^\op}\to
\calH^n_\C$ are inverse dualities, which preserve length, and in
particular give rise to a bijection between the simple $\C$-modules
and simple $\C^\op$-modules of maximal projective dimension $n$.
\end{enumerate}
\end{cor}
We end this section with an example of a generalized Artin-Schelter
regular category pointed out to us by Osamu Iyama.
\begin{example}\label{exam:Osamu}
  Let $\L$ be a non-semisimple finite dimensional algebra over a field
  $k$. Then a finitely generated left $\L$-module $M$ is called an
  \emph{$n$-cluster tilting module} ($(n-1)$-maximal orthogonal
  module, see \cite{I5,I6}) if
\begin{align}
\add M & = \{X\in\mod\L\mid \Ext^i_\L(X,M)=(0) \text{\ for\ }
  0<i<n\}\notag\\
 & = \{X\in\mod\L\mid \Ext^i_\L(M,X)=(0) \text{\ for\ } 0<i<n\}.\notag
\end{align}
Then $\add M$ is clearly a Krull-Schmidt category. Since $\Mod(\add
M)$ is equivalent to $\Mod\End_\L(M)$, the category $\add M$ is
Noetherian. By \cite[Proposition 3.5.1]{I5} it follows that $\add M$
is a generalized Artin-Schelter regular category of dimension $n+1$.
\end{example}

\section{Properties of generalized Artin-Schelter regular
  categories}\label{section:6}
This section is the categorical version of Section \ref{section1.5}
discussing Serre duality formulas and the relationship with finite
length Frobenius categories. We show that analogous results as for the
algebra setting are true in categorical setting, in particular that
the $\Ext$-category of any set of simple $\C$-modules in $\calH^n_\C$
for a generalized Artin-Schelter regular $K$-category of dimension $n$
with duality is a finite length Frobenius $K$-category. In addition,
in the Koszul case we give a converse of this.

First we consider the Serre duality formulas for generalized
Artin-Schelter regular categories. Let $\C$ be a generalized
Artin-Schelter regular $K$-category of dimension $n$. Then we have
seen that the subcategories $\calH^n_\C$ and $\calH^n_{\C^\op}$
consist exactly of the objects of finite length with composition
factors consisting of only simple $\C$-modules of maximal projective
dimension $n$. Moreover, $\calH^n_\C$ and $\calH^n_{\C^\op}$ are
abelian categories and $\tr^n_\C\colon \calH^n_\C\to \calH^n_{\C^\op}$
is an exact duality. If in addition all simple $\C$-modules $S$ are
finite dimensional over $K$ (in particular $S\simeq D^2(S)$), the
duality sends a simple $\C$-module to a simple $\C^\op$-module, but
the duality $D$ does not necessarily induce a duality from
$\calH^n_\C$ to $\calH^n_{\C^\op}$. In any case we then say that $\C$
is a generalized Artin-Schelter regular $K$-category of dimension $n$
\emph{with duality $D$}. As for algebras we have the following Serre
duality formulas.

\begin{lem}\label{lem:Serreduality}
Let $\C$ be a generalized Artin-Schelter regular $K$-category of dimension $n$
with duality $D$.
\begin{enumerate}
\item[(a)] The category $\C^\op$ is a generalized Artin-Schelter regular
 $K$-category of dimension $n$.
\item[(b)] For all $i$ there are natural isomorphisms
\[\Ext^i_{\Mod(\C)}(M,L)\simeq \Tor^{\Mod(\C)}_{n-i}(\tr^n_\C(M),L)\]
and
\[\varphi_i\colon D(\Ext^i_{\Mod(\C)}(M,L))\to
\Ext^{n-i}_{\Mod(\C)}(L,D\tr^n_\C(M))\]
for all $i\geq 0$, for any $L$ in $\Mod(\C)$ and any object $M$ in
$\calH^n_\C$.
\end{enumerate}
\end{lem}
\begin{proof} (a) From the definition it is clear that we only need to
show that the category of $\C^\op$-modules has global dimension
$n$. For the graded case this is a consequence of \cite[Theorem
1.18]{MVS1}. For the other case it follows in a similar
way.

(b) The proof is literary the same as in the algebra case.
\end{proof}

Let $\C$ be a generalized Artin-Schelter regular $K$-category of
dimension $n$ with duality $D$. Denote by $\T$ a full additive
subcategory generated by some simple $\C$-modules of maximal
projective dimension $n$, and consider the associated $\Ext$-category
$E(\T)$ of $\T$. Recall that $E(\T)$ has the same objects as $\T$,
while the morphisms for $A$ and $B$ in $E(\T)$ are given by 
\[\Hom_{E(\T)}(A,B)=\oplus_{i\geq 0}\Ext^i_{\Mod(\C)}(A,B).\]
Then define the natural functors
\[F\colon \Mod(\C)\to \Gr(E(\T))\]
and
\[F'\colon \Mod(\C)\to \Gr(E(\T)^\op)\]
given by
\[F(G)=\Ext^*_{\Mod(\C)}(G,-)\colon E(\T)\to \Gr(K)\]
and
\[F'(G')=\Ext^*_{\Mod(\C)}(-,G')\colon E(\T)^\op\to \Gr(K).\]
Then we can consider the following diagram of functors
\[\xymatrix@C=60pt{
\Mod(\C) \ar[r]^{F}\ar[d]_{D\tr^n_\C} &
\Gr(E(\T))\ar[d]^D\\
\Mod(\C)\ar[r]^{F'} & \Gr(E(\T^\op))
}\]
The following result shows that under certain stability conditions we
obtain a Frobenius category from $\T$.
\begin{thm}\label{thm:commdiagfunctors}
Let $\C$ be a generalized Artin-Schelter regular $K$-category of dimension $n$
with duality $D$. Let $\T$ be a full additive subcategory generated by
some simple $\C$-modules of maximal projective dimension $n$ such that
$D\tr^n_\L$ maps $\T$ into itself. Denote by $\F(\T)$ the full
subcategory of $\calH^n_\C$ consisting of objects with finite length
and compositions factors only in $\T$. 
\begin{enumerate}
\item[(a)] There is a commutative diagram of functors
\[\xymatrix@C=60pt{
\F(\T) \ar[r]^{F}\ar[d]_{D\tr^n_\C} &
\gr(E(\T))\ar[d]^D\\
\F(\T)\ar[r]^{F'} & \gr(E(\T^\op))
}\]
where $\gr E(\T)$ is a finite length category.
\item[(b)] Let $\T$ be such that $D\tr^n_\C$ induces a bijection on
the simple $\C$-modules in $\T$. Then the graded $K$-category
$\gr(E(\T))$ is Frobenius.
\item[(c)] Let $\T$ be such that $D\tr^n_\C$ induces a bijection on
the simple $\C$-modules in $\T$. There is a commutative diagram of
functors
\[\xymatrix@C=60pt{
\F(\T)\ar[r]^{F}\ar[d]_{D\tr^n_\C} &
\gr(E(\T))\ar[d]^{(-)^* D}\\
\F(\T)\ar[r]^{F} & \gr(E(\T))
}\]
\end{enumerate}
\end{thm}
\begin{proof} (a) For a projective functor $\Ext^*_{\Mod(\C)}(S,-)$ in $\gr
E(\T)$ we have
\[\Ext^*_{\Mod(\C)}(S,-)=\amalg_{i=0}^n\Ext^i_{\Mod(\C)}(S,-),\]
where
\[\Ext^i_{\Mod(\C)}(S,-)\simeq \Hom_{\Mod(\C)}(P_i/\rad P_i,-)\]
when
\[0\to P_n\to P_{n-1}\to \cdots \to P_1\to P_0\to S\to 0\]
is a projective resolution of $S$. Since $P_i$ is finitely generated,
the support of the functor $\Hom_{\Mod(\C)}(P_i/\rad P_i,-)$ is
finite. Hence the support of $\Ext^*_{\Mod(\C)}(S,-)$ is finite and
therefore $\gr E(\T)$ is a finite length category.

Since the functor $F$ maps semisimple $\C$-modules in $\F(\T)$ to
finitely generated projective objects in $\gr E(\T)$, $F$ is a half
exact functor and $\gr E(\T)$ is a finite length category, therefore 
$F|_{\F(\T)}$ has its image in $\gr E(\T)$. Similar arguments are
used for $F'$. Hence the functors in the diagram end up in the
categories indicated.

By Lemma \ref{lem:Serreduality} the image of the functors $DF$ and
$F'D\tr^n_\L$ on a $\C$-module $M$ in $\F(\T)$ are isomorphic as
graded vector spaces. So we need to show that they in fact are
isomorphic as $E(\T^\op)$-modules. Recall that the isomorphism from
$DF(M)$ to $(F'D\tr^n_\L)(M)$ as graded vector spaces is induced by
the isomorphisms
\[\varphi_i\colon  D(\Ext^{n-i}_{\Mod(\C)}(M,L))\to
\Ext^i_{\Mod(\C)}(L,D\tr^n_\C(M))\]
for all $i\geq 0$. For any element $\theta$ in $\Ext^j_{\Mod(\C)}(T,T)$
we want to prove that $\varphi_{i+j}(f)\theta=\varphi_i(f\theta)$ for
all $f$ in $D\Ext^{i+j}_{\Mod(\C)}(M,T)$.

Represent $\theta$ in $\Ext^j_{\Mod(\C)}(T,T)$ as a homomorphism
$\theta\colon \Omega^j_\C(T)\to T$. This gives rise to the
following diagram
\begin{small}
\[\xymatrix{
D\Ext^{i+j}_{\Mod(\C)}(M,T)
\ar[r]^{D\Ext^{i+j}_{\Mod(\C)}(M,\theta)}
\ar[d]^{\varphi_{i+j}} &
D\Ext^{i+j}_{\Mod(\C)}(M,\Omega^j_\C(T)) \ar[r] \ar[d]^{\varphi_{i+j}} &
D\Ext^i_{\Mod(\C)}(M,T)\ar[d]^{\varphi_i}\\
\Ext^{n-i-j}_{\Mod(\C)}(T,D\tr^n_\C(M))
\ar[r]^{\Ext^{n-i-j}_{\Mod(\C)}(\theta,D\tr^n_\C(M))} &
\Ext^{n-i-j}_{\Mod(\C)}(\Omega^j_\C(T),D\tr^n_\C(M)) \ar[r] &
\Ext^{n-i}_{\Mod(\C)}(T,D\tr^n_\C(M))}\]
\end{small}%
The first square commutes, since $\varphi_{l}$ is a natural
transformation for all $l$. The horizontal morphisms in the second
square are compositions of connecting homomorphisms or dual
thereof. The isomorphisms $\varphi_l$ are all induced by the
isomorphism for $\Hom_\C(-,C)$ and any $\C$-module $G$
\[\psi\colon D\Hom_{\Mod(\C)}({_\C(-,C)},G)\to
\Hom_{\Mod(\C)}(G,D({_\C(-,C)}^*))\]
given by $\psi(g)_X(u)(s)=g(\alpha_X(u)\comp {_\C(-,s)})$ for $g$ in
$D(\Hom_{\Mod(\C)}({_{\Mod(\C)}(-,C)},G))$, $u$ in $G(X)$ and $s$ in
$\Hom_\C(C,X)$ for all $X$ in $\C$, where
\[\alpha_X\colon  G(X)\to \Hom_{\Mod(\C)}({_\C(-,X)},G)\]
is the Yoneda isomorphism. This is a natural isomorphism in both
variables. It follows from this that the isomorphisms $\varphi_l$ all
commute with connecting morphisms, that is, if $0\to X\to Y\to Z\to 0$
is exact, then the diagram
\[\xymatrix{
D\Ext^{i+1}_{\Mod(\C)}(M,X)\ar[r]^{D(\partial)} \ar[d]^{\varphi_{i+1}} &
D\Ext^i_{\Mod(\C)}(M,Z)\ar[d]^{\varphi_i}\\
\Ext^{n-i-1}_{\Mod(\C)}(X,D\tr^n_\C(M))\ar[r]^{\partial'} &
\Ext^{n-i}_{\Mod(\C)}(Z,D\tr^n_\C(M))}\]
is commutative.

Then going back to our first diagram, which by the above is
commutative, the upper diagonal path is computing $\varphi_i(f\theta)$
while the lower diagonal path is computing $\varphi_{i+j}(f)\theta$
both for any $f$ in $D\Ext^{i+j}_{\Mod(\C)}(M,T)$. Hence we conclude that
they are equal and $DF(M)$ and $F'(D\tr^n_\C(M))$ are isomorphic as
$E(\T^\op)$-modules for all $M$ in $\F(\T)$.

(b) By (a) we have that
\[D(\Ext^*_{\Mod(\C)}(-,D(\tr^n_\C(S)))) \simeq
\Ext^*_{\Mod(\C)}(S,-)\]
as objects in $\gr(E(\T))$. Since $D\tr^n_\C$ is a bijection on the
simple objects in $\T$, it follows that the projective and the
injective objects in $\gr E(\T)$ coincide, and $\gr E(\T)$ is
Frobenius.

(c) Recall that $F'(G')=\Ext^*_{\Mod(\C)}(-,G')\colon E(\T^\op)\to
\Gr(K)$ for all $G'$ in $\F(\T)$, and that $(F'(G'))^*$ is given
by
\[(F'(G'))^*(S)=\Hom_{E(\T)^\op}(F'(G'),\Ext^*_{\Mod(\C)}(-,S))\]
for $S$ in $E(\T)$. For $G'=S'$ in $\T$ we have that
\begin{align}
(F'(S'))^*(S) &
 =\Hom_{E(T)^\op}(\Ext^*_{\Mod(\C)}(-,S'),\Ext^*_{\Mod(\C)}(-,S))\notag\\
 & \simeq \Ext^*_{\Mod(\C)}(S',S)\notag\\
 & = F(S')(S)\notag
\end{align}
Then using similar arguments as for the algebra situation, one can
show that $(-)^*F'\simeq F\colon \F(\T)\to \gr E(\T)$. The claim
in (c) follows from this.
\end{proof}
The above result has the following immediate consequence.
\begin{cor}
Let $\L$ be a non-semisimple finite dimensional algebra over a field
$k$, and let $M$ be an $n$-cluster tilting $\L$-module. Denote by $\G$
the endomorphism ring $\End_\L(M)$ of $M$. Let $T=\G/\rrad_\G$. Then
$\Ext^*_\G(T,T)$ is a finite dimensional Frobenius algebra.
\end{cor}
Theorem \ref{thm:commdiagfunctors} shows that for any generalized
Artin-Schelter regular $K$-category $\C$ with duality and any full
subcategory $\T$ of simple $\C$-modules of maximal projective
dimension $n$ such that $D\tr^n_\C$ induces a bijection on the simple
$\C$-modules in $\T$, the associated $\Ext$-category $E(\T)$ is
Frobenius. Here $\C$ need not to be Koszul. In the Koszul case the
converse is also true as we prove next.

\begin{thm}
  Let $\C$ be a Koszul $K$-category, such that $\gr E(\S (\C))$ is
  indecomposable for the $\Ext$-category $E(\S(\C))$ of the full 
  subcategory $\S(\C)$ of $\C$-modules generated by the simple
  objects. Assume that $\gr E(\S (\C))$ is Frobenius and a finite
  length category with Loewy length $n$. Then $\C$ is generalized
  Artin-Schelter regular of dimension $n$ and all simple objects have
  the same projective dimension.
\end{thm}
\begin{proof}
The proof will be the same as in the algebra case given in \cite{M},
it will rely on the use of Koszul duality. We give it here for
completeness.

Let $\D$ be an indecomposable Frobenius Koszul $K$-category which is a
finite length category. First we prove that all indecomposable
projective objects in $\gr(\D)$ have the same Loewy length. Let
$P=\Hom_{\D}(-,C)$ for an indecomposable object $C$ in $\D$ of maximal
Loewy length. Let $Q=\Hom_{\D}(-,B)$ be an indecomposable projective
in $\gr(\D)$ not isomorphic to $P$ with $Q/\rad Q$ as a direct summand
of $\Hom_{\D}(-,C)_{1}$. By Yoneda's Lemma there is a morphism
$\widetilde{f}\colon Q\rightarrow \rad P$ induced by a morphism
$f\colon B\rightarrow C$ in $\Hom_{\D}(B,C)_{1}$. Since $P$ and $Q$
both are indecomposable projective injective objects, we infer that
the socle of $Q$ is contained in $\Ker\widetilde{f}$. It follows that
$\LL (Q)>\LL(\rad P)=\LL(P)-1$, where $\LL(H)$ denotes the Loewy
length of a functor $H$. Since $\LL(P)$ is maximal, $\LL(Q)=\LL(P)$.

Let $Q=\Hom_\D(-,B)$ be a indecomposable summand of the injective
envelope of $P/\soc P$. As above there is a morphism $f\colon C\to B$
in $\Hom_\D (C,B)_1$. Similarly as above using the maximality of
$\LL(P)$ we conclude that $\LL(Q)=\LL(P)$. We obtain by induction that
for all indecomposable objects $B$, which has a chain of non-zero maps
(in either directions) to $C$, the projective objects $\Hom_\D(-,B)$
has the same Loewy length as $\Hom_\D (-,C)$. Since $\D$ is
indecomposable, all (indecomposable) projective objects have the same
Loewy length.

Now, since $\D=\E(\S (\C))$ is assumed to be Frobenius and $\C$ is a
Koszul $K$-category, it follows by the proof of \cite[Theorem 2.3 (e),
(f)]{MVS1} that all the simples in $\C$ have projective dimension
$n$. By \cite[Theorem 1.18] {MVS1} the category $\C$ has global
dimension $n$.

Let $C$ be an indecomposable object in $\C$, and let $S_{C}=\Hom_{\C}(-,C)/%
\rad(-,C)$. We want to show that $\Ext_{\Mod(\C)}^{i}(S_{C},\Hom_{\C%
}(-,X))=(0) $ for $i\neq n$ for all objects $X$ in $\C$. Let
\begin{equation*}
0\rightarrow P_{n}[-n]\rightarrow P_{n-1}[-n+1]\rightarrow \cdots
\rightarrow P_{1}[-1]\rightarrow P_{0}\rightarrow S_{C}\rightarrow 0
\end{equation*}
be a minimal graded projective resolution of $S_{C}$.

Let $Q=\Hom_{\C}(-,X)$ be an indecomposable projective in $\gr(\C)$. Recall
that we considered the functor $E\colon \gr(\C)\rightarrow \Gr E(\S (\C))$
given by
\[E(H)=\Ext_{\Mod(\C)}^{\ast }(H,-)\colon E(\S (\C))\rightarrow \Gr(K)\]
induces a duality between the subcategory of linear functors
$\mathcal{L}(\C)$ in $\Mod(\C)$, and the subcategory of linear
functors $\mathcal{\D}$ in $\Mod(\D)$.

An element of $\Ext_{\Mod(\C)}^{i}(S_{C},\Hom_{\C}(-,X))$ is represented
by a degree zero map $h\colon\Omega^{i}(S_{c})\rightarrow Q[-m]$,
where $Q=\Hom_{\C}(-,X)$. Since $E(\rad^{i}Q[i])\simeq
\Omega^{i}E(Q)[i]$ and $E(Q)$ is simple, $E(\rad^{i}Q)$ is
indecomposable and therefore $\rad^{i}Q$ is indecomposable.
Similarly, $E(\Omega^{i}(S_{c})[i])\simeq \rad^{i}E(S_{C})[i]$ is a
subfunctor of an indecomposable projective injective functor of finite
length, hence with simple socle. It follows that any subfunctor of
$E(S_{C})$ is indecomposable, in particular, $\rad^{i}E(S_{C})$ is
indecomposable, then $\Omega^{i}(S_{c}) $ is indecomposable.
From these observations and the fact that projective dimension of
$S_{C}$ is $n$, it follows the map $h$ is not an epimorphism, hence
$m<i$ and $\Im h\subset \rad^{i-1-m}Q[-m]$. Consider the following push
out
\[\xymatrix{
0\ar[r] & \Omega ^{i}(S_{c}) \ar[r]\ar[d]^h
        & P_{i-1}[-i+1] \ar[r]\ar[d]
        & \Omega ^{i-1}(S_{c}) \ar[r]\ar@{=}[d] & 0\\
0\ar[r] & \rad^{i-1-m}Q[-m]\ar[r]
        & Z\ar[r]
        & \Omega^{i-1}(S_{c})\ar[r] & 0}\]
Denote the lower exact sequence in the above diagram by $\theta$.
The end terms in $\theta$ are linear functors generated in
degree $i-1$, therefore $Z$ is a linear functor generated in degree
$i-1$. Hence the exact sequence, $\theta[i-1]$,
\[0\to \rad^{i-1-m}Q[i-1-m]\to Z[i-1]\to \Omega^{i-1}(S_{c})[i-1]\to 0\]
is a sequence of linear functors. Applying Koszul duality we obtain an
exact sequence linear functors
\[0\to E(\Omega^{i-1}(S_{c})[i-1])\to E(Z[i-1])\to
E(\rad^{i-1-m}Q[i-1-m])\to 0\]
in $\Mod(\D)$, which is isomorphic to the following sequence
\[\theta'\colon 0\to \rad^{i-1}E(S_{c})[i-1])\extto{g}
E(Z[i-1])\to E(\rad^{i-1-m}Q[i-1-m])\to 0\]
There is an inclusion map $j\colon \rad^{i-1}E(S_{c})[i-1])\to
E(S_{c})[i-1]$ and $E(S_{c})[i-1]$ it is projective injective, it
follows there exists a map $f\colon E(Z[i-1]\to E(S_{c})[i-1]$ such
that $fg=j$ the image of f is contained in
$\rad^{i-1}E((S_{c}))[i-1]$. Therefore the exact sequence $\theta'$
splits. By Koszul duality the exact sequence $\theta$ also splits.
We have proved the map $h\colon \Omega ^{i}(S_{c})\to Q[-m]$ extends
to $P_{i-1}[-i+1]$, and the corresponding extension splits.

If $\Hom(S_{C},Q)\neq (0)$, then there is an integer $k$ such that
$S_{C}$ is a summand of $\rad^{k}Q$, but we saw above that $\rad^{k}Q$
is indecomposable. This implies $\rad^{k}Q=S$ and $Q$ is of finite
length, which would imply $E(Q)$ has finite projective dimension
contradicting $\D$ is Frobenius.

Dualizing with $(-)^{\ast }$ the exact sequence
\[0\to P_{n}[-n]\to P_{n-1}[-n+1]\to \cdots \to P_{1}[-1]\to P_{0}\to
S_{C}\to 0\]
and shifting by $n$, we obtain the exact sequence
\[0\to P_{0}^{\ast }[-n]\to P_{1}^{\ast }[-n+1]\to
\cdots \to P_{n-1}^{\ast }[-1]\to P_{n}^{\ast }\to
M\to 0\]
Hence, since $P_{n}$ is indecomposable, the functor $M$ is linear of
projective dimension $n$ and indecomposable.

Since the opposite category is Koszul, we can apply Koszul duality to
$M$ to get an indecomposable functor $E(M)$ in $\Gr(\D)$ of Loewy
length $n$, but the only indecomposable functors of length $n$ in
$\Gr(\D^{\op}$ are the projective.  It follows $M$ is simple.
We have proved $\Ext_{\Mod(\C)}^{n}(S_{C},\Hom_{\C}(-,X))$ is simple.
\end{proof}

Let us now briefly return to the polynomial ring $\L=K[x_1,x_2,\ldots,
x_n]$ in $n$ indeterminants $x_1,x_2,\ldots,x_n$ over an algebraically
closed field $K$ discussed in Section \ref{section2}.  Proposition
\ref{prop:commdiag} (b) in particular says that the Yoneda algebra
$\amalg_{i\geq 0}\Ext^i_\L(K,K)$ is Frobenius. This is well-known and
it is in fact isomorphic to the exterior algebra on $K^n$.  However we
can consider $\L$ as an ungraded $K$-algebra, and we have seen that
$\L$ is a generalized Artin-Schelter regular $K$-algebra. If $\C$ is
the category with one object $*$ with endomorphism ring $\L$, then the
category of simple objects $E(\S)$ in $\Mod(\C)$ is the category of
simple $\L$-modules $S_{\overline{a}}=\L/\frakm_{\overline{a}}$ with
$\frakm_{\overline{a}}=(x_1-a_1,x_2-a_2,\ldots,x_n-a_n)$ and
homomorphisms
\[\Hom_{E(\S)}(S_{\overline{a}}, S_{\overline{b}})=
\amalg_{i\geq 0} \Ext^i_\L(S_{\overline{a}},S_{\overline{b}}).\]
Each $\Ext^i_\L(S_{\overline{a}}, S_{\overline{b}})$ is a module over
$\L$, where the action of $\L$ can be taken through $S_{\overline{a}}$
or $S_{\overline{b}}$. Hence
$\Ext^i_\L(S_{\overline{a}},S_{\overline{b}})$ is annihilated by
$\frakm_{\overline{a}}+\frakm_{\overline{b}}=\L$ if $\overline{a}\neq
\overline{b}$. Therefore $\Gr(E(\S))$ splits in a product of
categories with each one isomorphic to the category of modules over
$\amalg_{i\geq 0}\Ext^i_\L(S_{\overline{a}},S_{\overline{a}})$. Since
$S_{\overline{a}} \simeq K^{\sigma_{\overline{a}}}$, it follows that
each of these algebras are isomorphic to the exterior algebra.

\section{Generalized Artin-Schelter regular categories of global
  dimension $2$}\label{section:7}
\sloppy The main example providing the motivation for this work is the
category $\Mod(\mod\L)$ of all additive functors from
$(\mod\L)^\op$ to $\Mod K$ for a finite dimensional $K$-algebra
$\L$. However, in the our definition of a generalized Artin-Schelter
regular $K$-category $\C$, we require that the category is Noetherian,
that is, $\Mod\C$ is Noetherian. It is shown by Auslander in
\cite[Theorem 3.12]{A2} that $\Mod(\mod\L)$ is Noetherian if and only
if $\L$ is of finite representation type. However, if we consider the
associated graded category, $\A_\gr(\mod\L)$, of $\mod\L$, that is,
the objects are the same as for $\mod\L$ and the morphism is given by
\[\Hom_{\A_\gr(\mod\L)}(A,B)=\oplus_{i\geq 0}
\rad^i_\L(A,B)/\rad^{i+1}_\L(A,B).\] Then consider the graded
$K$-category $\Gr(\A_\gr(\mod\L))$ of graded functors from
$(\A_\gr(\mod\L)^\op\to \Gr(K)$ (see \cite{MVS1}). Then by \cite{IT}
the category $\C=\Gr(\A_\gr(\mod\L))$ is a generalized Artin-Schelter
regular $K$-category of dimension $2$. Furthermore,
$\Gr(\A_\gr(\mod\L))$ gives rise to finite length Frobenius categories
through naturally associated $\Ext$-categories. This is why we discuss
Krull-Schmidt categories $\C$ where all finitely generated projective
objects in $\Mod(\C)$ has finite length and generalized Artin-Schelter
regular $K$-categories of dimension $2$ in this section. We show that
an indecomposable positively graded locally finite Krull-Schmidt
$K$-category generated in degrees $0$ and $1$ with finite presented
simple objects and finitely generated projective objects in $\Mod(\C)$
have finite length, then the number of isomorphism classes of
indecomposable objects in $\C$ is countable. In addition, a
generalized Artin-Schelter regular positively graded locally finite
$K$-category generated in degrees $0$ and $1$ of dimension $2$ is
coherent and $\gr(\C)$ is abelian.

We start by showing that the number of isomorphism classes of
indecomposable objects in a positively graded locally finite
$K$-category generated in degrees $0$ and $1$ with finitely presented
simple objects in $\Mod(\C)$ is countable.

\begin{lem}\label{lem:countable}
Let $\C$ be a positively graded locally finite $K$-category generated
in degrees $0$ and $1$ such that the graded simple objects
$S_C=\Hom_\C(-,C)_0$ are finitely presented. Assume that $\C$ is
indecomposable and that each projective object has finite length. Then
the isomorphism classes of indecomposable objects in $\C$ is
countable.
\end{lem}
\begin{proof}
Let $C$ be an object in $\C$ and denote by $\langle C\rangle$ the full
subcategory of $\C$ consisting of all objects $X$ such that there is a
chain of morphisms in $\C$
\[\xymatrix{X=X_0\ar@{-}[r]^-{f_1} & X_1\ar@{-}[r]^{f_2} &
  X_2\ar@{..}[r] & X_{n-1}\ar@{-}[r]^-{f_{n-1}} & X_n=C},\]
where $f_i$ is either in $\Hom_\C(X_{i-1},X_i)_1$ or
$\Hom_\C(X_i,X_{i-1})_1$.

The fact that simple objects are finitely presented means that there
are only a finite number of objects $X$ connected to an object $Y$
with a map in degree $1$.

Since each projective object has finite length, and $\C$ is generated
in degrees $0$ and $1$, then each homogeneous morphism $f\colon X\to
C$ in $\C$ with a positive degree $l$ is a sum of composition of maps
$f_1f_2\cdots f_l$ with $f_i$ in degree $1$. Then $\langle C\rangle$
is countable. If $D$ is in $\C\setminus\langle C\rangle$, then by the
definition of $\langle C\rangle$ there is no morphism from $D$ to any
object of $\langle C\rangle$ and from any object of $\langle C\rangle$
to $D$. This contradicts the fact that $\C$ is connected, so that we
have $\C=\langle C\rangle$ and $\C$ is countable.
\end{proof}
Note that when $\C$ is a generalized Artin-Schelter regular Koszul
$K$-category of dimension $n$ with duality $D$ and $\T$ is a set of
simple $\C$-modules in $\calH^n_\C$ permuted by $D\tr^n_\C$, then the
$\Ext$-category $E(\T)$ satisfies the assumptions in Lemma
\ref{lem:countable}.

We end this section and this paper with showing that generalized
Artin-Schelter regular $K$-categories of dimension $2$ generated in
degrees $0$ and $1$ are coherent. To this end recall the definition of
the functor $t^n\colon \Mod(\C)\to \Mod(\C)$ preceding Lemma
\ref{lem:torsionradical}. Since we are discussing the global dimension
$2$ case in this section, we are considering the functor $t^2$.

\begin{thm}
  Let $\mathcal{C}$ be a generalized Artin-Schelter regular positively
  graded locally finite $K$-category generated in degrees $0$ and $1$
  of global dimension $2$. Then the following statements are true.
\begin{enumerate}
\item[(a)] If $F$ is a finitely generated functor, then $t^2(F)$ is of
  finite length.
\item[(b)] The category $\C$ is coherent.
\end{enumerate}
\end{thm}
\begin{proof}
It was proved in Lemma \ref{lem:torsionradical} that $t^{n}$ is a
radical. Furthermore, for any (finitely presented) functor $F$, the
functor $F/t^{2}(F)$ is $\calH^2_\C$-torsion free and therefore has a
projective dimension at most $1$ by Lemma \ref{lem:projdimtorfree}.

Consider the commutative exact diagram
\[\xymatrix{
          &             & 0\ar[d] & 0\ar[d] & & \\
& 0\ar[r] & \Omega F \ar[r]^j\ar[d]
                       & \Omega (F/t^{2}F)\ar[r]\ar[d]
                       & t^{2}F\ar[r] & 0 \\
          &             & P\ar[d]\ar@{=}[r] & P\ar[d] & & \\
0\ar[r] & t^2(F)\ar[r] & F\ar[r]\ar[d]
                       & F/t^2(F)\ar[r]\ar[d]
                       & 0 & \\
        &              & 0 & 0 & &
}\]
where $\Omega (F/t^{2}F)\simeq \oplus_{i\in I} Q_{i}$ and each $Q_{i}$
projective. Since $\Omega F$ is finitely generated, there exists a
finite subset $J\subset I$ such that $j(\Omega F)\subset
\oplus_{i\in J}Q_{i}$.

We have the following commutative exact diagram
\[\xymatrix{
        & & 0\ar[d] & 0\ar[d] & \\
0\ar[r] & \Omega F\ar[r]\ar@{=}[d]
        & \oplus_{i\in J} Q_i\ar[r]\ar[d]
        & Z\ar[r]\ar[d] & 0\\
0\ar[r] & \Omega F\ar[r]
        & \oplus_{i\in I} Q_i\ar[r]\ar[d]
        & t^2(F)\ar[r]\ar[d] & 0\\
        & & \oplus_{I\setminus J} Q_i\ar[d]\ar@{=}[r]
          & \oplus_{I\setminus J} Q_i\ar[d] & \\
        & & 0 & 0 & \\
}\]
each $Q_{i}$ with $i$ in $I\setminus J$ is projective and
$\calH^2_\C$-torsion, therefore it has a simple subobject $S$ from
$\calH^2_\C$ and $L=Q_{i}/S$ has projective dimension $3$, a
contradiction. It follows that $I$ is finite and that $t^{2}F$ has
finite length. By the Horseshoe Lemma, $F$ has a projective resolution
consisting of finitely generated projective. It follows $\mathcal{C}$
is coherent.
\end{proof}
This has the following immediate consequence.
\begin{cor}
Let $\mathcal{C}$ be a generalized Artin-Schelter regular positively graded
locally finite $K$-category generated in degrees $0$ and $1$ of global
dimension $2$. Then the category of finitely presented functors in
$\Gr(\C)$ is abelian.
\end{cor}


\begin{thebibliography}{99}
\bibitem{AF} Anderson, F., Fuller, K., \emph{Rings and categories of
    modules}, Graduate Texts in Mathematics, Vol.\
  13. Springer-Verlag, New-York-Heidelberg, 1992.
\bibitem{ArS} Artin, M., Schelter, W., \emph{Graded algebras of global
    dimension $3$}, Adv.\ Math., 66 (1987), 171--216.
\bibitem{A2} Auslander, M., \emph{A functorial approach to
    representation theory}, Representations of algebras (Puebla,
  1980), 105--179, Lecture Notes in Math., 944.
\bibitem{A} Auslander, M., \emph{Representation theory of artin
    algebras \emph{I}}, Comm.\ Algebra 1 (1974), 177--268.
\bibitem{AR} Auslander, M., Reiten, I., \emph{Representation theory of
    artin algebras \emph{III}. Almost split sequences}, Communications in
  Algebra \textbf{3} (1975), 239-294.
\bibitem{ARS} Auslander, M., Reiten, I., Smal\o , S.\
  O., \emph{Representation Theory of Artin Algebras}, Cambridge
  Studies in Advanced Mathematics, \textbf{36}, Cambridge University
  Press, Cambridge, (1995).
\bibitem{BGL} Baer, D., Geigle, W., Lenzing, H., \emph{The
    preprojective algebra of a tame hereditary Artin algebra}, Comm.\
  Alg.\ 15 (1987), 425--457.
\bibitem{GMT} Guo, J.\ Y., Martin\'ez-Villa, R., Takane, M.,
  \emph{Koszul generalized Auslander regular algebras}, Algebras and
  modules, II (Geiranger, 1996), 263--283, CMS Conf.\ Proc., 24,
  Amer.\  Math.\ Soc., Providence, RI, 1998.
\bibitem{IT} Igusa, K., Todorov, G., \emph{Radical layers of
  representable functors}, J.\ Algebra 89 (1984), no.\ 1, 105--147.
\bibitem{I5} Iyama, O., \emph{Higher-dimensional Auslander-Reiten
    theory on maximal orthogonal subcategories}, Adv.\ Math.\ 210
  (2007), no.\ 1, 22–-50.
\bibitem{I6} Iyama, O., \emph{Auslander correspondence}, Adv.\ Math.\
  210 (2007), no.\ 1, 51-–82.
\bibitem{I1} Iyama, O., \emph{$\tau$-categories \emph{I}: Ladders},
  Algebras and Repre.\ Theory,  (2005) 8, 297-–321.
\bibitem{I2} Iyama, O., \emph{$\tau$-categories \emph{II}: Nakayama Pairs
and Rejective Subcategories}, Algebras and Repre.\ Theory,
  (2005) 8, 449-–477.
\bibitem{I3} Iyama, O., \emph{$\tau$-categories \emph{III}: Auslander
    orders and Auslander-Reiten quivers}, Algebras and Repre.\ Theory,
  (2005) 8, 601--619.
\bibitem{I4} Iyama, O., \emph{Symmetry and duality on $n$-Gorenstein
    rings}, J.\ of Algebra 269 (2003) 528-–535.
\bibitem{M1} Martin\'ez-Villa, R., \emph{Graded, selfinjective, and
    Koszul algebras}, J. Algebra 215 (1999), no. 1, 34--72. 
\bibitem{M2} Martin\'ez-Villa, R., \emph{Applications of Koszul
    algebras: the preprojective algebra}, Representation theory of
  algebras (Cocoyoc, 1994), 487--504, CMS Conf.\ Proc., 18,
  Amer.\ Math.\ Soc., Providence, RI, 1996. 
\bibitem{M3} Martin\'ez-Villa, R., \emph{Skew group algebras and their
    Yoneda algebras}, Math.\ J.\ Okayama Univ.\ 43 (2001), 1--16. 
\bibitem{MV} Martin\'ez Villa, R., \emph{Koszul algebras and the
Gorenstein condition}, Representations of algebras (S\~ao Paulo, 1999),
135--156, Lecture Notes in Pure and Appl.\ Math., 224 (2002).
\bibitem{MVS1} Martin\'ez Villa, R., Solberg, \O., \emph{Graded and
Koszul categories}, preprint 2008.
\bibitem{M} Mitchel, \emph{Rings with several objects},
  Advances in Math.\ 8 (1972), 1--161.
\bibitem{LPWZ} Lu, D. M., Palmieri, J. H., Wu, Q. S., Zhang, J. J.,
\emph{$A_\infty$-algebras for ring theorists},
Proceedings of the International Conference on Algebra.
Algebra Colloq.\ 11 (2004), no.\ 1, 91--128.
\bibitem{RvdB} Reiten, I., van den Bergh, M., \emph{Noetherian
    hereditary abelian categories satisfying Serre duality}, J.\
  Amer.\ Math.\ Soc.\ 15 (2002), no.\ 2, 295--366. 
\bibitem{S} Smith, P., \emph{Some finite dimensional algebras related
to elliptic curves}, Representation theory of algebras and related
topics (Mexico City, 1994), 315--348, CMS Conf.\ Proc., 19,
Amer.\ Math.\ Soc., Providence, RI, 1996.
\end{thebibliography}
\end{document}